\documentclass[12pt]{amsart}

\newtheorem{theorem}{Theorem}[section]
\newtheorem{lemma}[theorem]{Lemma}
\newtheorem{proposition}[theorem]{Proposition}
\newtheorem{corollary}[theorem]{Corollary}

\theoremstyle{remark}
\newtheorem*{remark}{Remark}
\newtheorem{remarks}{Remark}

\theoremstyle{definition}
\newtheorem*{notation}{Additional notation}
\newtheorem*{outline}{Outline of the proof}

\newcommand{\et}{\quad\mbox{and}\quad}
\newcommand{\bA}{\mathbb{A}}
\newcommand{\bC}{\mathbb{C}}
\newcommand{\bP}{\mathbb{P}}
\newcommand{\bQ}{\mathbb{Q}}
\newcommand{\bR}{\mathbb{R}}
\newcommand{\bZ}{\mathbb{Z}}
\newcommand{\cC}{{\mathcal{C}}}

\newcommand{\cK}{{\mathcal{K}}}
\newcommand{\cM}{{\mathcal{M}}}
\newcommand{\cO}{{\mathcal{O}}}
\newcommand{\cS}{{\mathcal{S}}}
\newcommand{\Cbar}{{\bar{\cC}}}
\newcommand{\Cg}{\cC^g}
\newcommand{\Cphi}{\cC^\varphi}
\newcommand{\dist}{\mathrm{dist}}
\newcommand{\fp}{{\mathfrak{p}}}
\newcommand{\fP}{{\mathfrak{P}}}
\newcommand{\fq}{{\mathfrak{q}}}

\newcommand{\GL}{\mathrm{GL}}

\newcommand{\Ibar}{{\bar{I}}}
\newcommand{\Kbar}{{\bar{K}}}
\newcommand{\Ktheta}{\cK^\theta}
\newcommand{\lambdabar}{{\bar{\lambda}}}
\newcommand{\lambdag}{\lambda^g}
\newcommand{\lambdaphi}{\lambda^\varphi}
\newcommand{\linpol}{A}

\newcommand{\rank}{\mathrm{rank}}
\newcommand{\rmH}{\mathrm{H}}

\newcommand{\uu}{\mathbf{u}}

\newcommand{\ux}{\mathbf{x}}
\newcommand{\uy}{\mathbf{y}}
\newcommand{\uz}{\mathbf{z}}
\newcommand{\uxi}{\underline{\xi}}

\newcommand{\vinM}{{v\in\cM}}
\newcommand{\vol}{\mathrm{Vol}}
\newcommand{\Vbar}{\overline{V}}

\begin{document}

\baselineskip=17pt

\title[Simultaneous approximation in transcendence degree one]
{Simultaneous approximation by conjugate algebraic numbers in
fields of transcendence degree one}
\author{Damien ROY}
\address{
   D\'epartement de Math\'ematiques\\
   Universit\'e d'Ottawa\\
   585 King Edward\\
   Ottawa, Ontario K1N 6N5, Canada}
\email{droy@uottawa.ca}
%
\subjclass{Primary 11J13; Secondary 11J82}
\thanks{Work partially supported by NSERC and CICMA}
\keywords{simultaneous approximation, conjugate algebraic numbers,
polynomials, Gel'fond's criterion, heights}

\begin{abstract}
We present a general result of simultaneous approximation to several
transcendental real, complex or p-adic numbers $\xi_1,...,\xi_t$ by
conjugate algebraic numbers of bounded degree over $\bQ$, provided
that the given transcendental numbers $\xi_1,...,\xi_t$ generate
over $\bQ$ a field of transcendence degree one.  We provide sharper
estimates for example when $\xi_1,...,\xi_t$ form an arithmetic
progression with non-zero algebraic difference, or a geometric
progression with non-zero algebraic ratio different from a root of
unity. In this case, we also obtain by duality a version of
Gel'fond's transcendence criterion expressed in terms of polynomials
of bounded degree taking small values at $\xi_1,...,\xi_t$.
\end{abstract}

\maketitle

\section{Introduction}
 \label{intro}
The basic problem of approximation to real numbers by algebraic
numbers of bounded degree has attracted much attention since the
pioneer work \cite{Wi} of E.~Wirsing in 1961.  In their seminal
paper \cite{DS} of 1969, H.~Davenport and W.~M.~Schmidt proposed
an innovative approach based on geometry of numbers which allowed
them to deal with approximation by algebraic integers. Recently,
Y.~Bugeaud and O.~Teuli\'e observed that it can also be used to
treat approximation by algebraic integers of a fixed degree
\cite{BT}. The sharpest result in this direction is due to
M.~Laurent \cite[Cor.]{La}.  Simplifying slightly, it shows that,
for each integer $n\ge 2$ and each real number $\xi$ which is not
algebraic over $\bQ$ of degree at most $n/2$, there are infinitely
many algebraic integers $\alpha$ of degree $n$ over $\bQ$ which
satisfy $|\xi-\alpha|\le \rmH(\alpha)^{-n/2}$, where the {\em
height} $\rmH(\alpha)$ of $\alpha$ is defined as the largest
absolute value of the coefficients of the irreducible polynomial
of $\alpha$ over $\bZ$. Similar estimates valid for a $p$-adic
number $\xi\in\bQ_p$ are also known \cite{Mo,Te}.

The present work deals with the problem of simultaneous
approximation to several numbers by conjugate algebraic numbers.
Naive heuristic arguments based on Dirichlet box principle suggest
that, for each integer $t\ge 1$ and each choice of transcendental
real numbers $\xi_1,\dots,\xi_t$, there exist constants $n_0\ge 1$
and $c>0$ with the property that, for each integer $n\ge n_0$,
there are infinitely many algebraic numbers $\alpha$ of degree $n$
over $\bQ$ which admit distinct real conjugates $\alpha_1, \dots,
\alpha_t$ with $|\xi_i-\alpha_i|\le \rmH(\alpha)^{-cn}$ for
$i=1,\dots,n$.  For $t=1$, this heuristic statement is true, by
the above, with $n_0=2$ and $c=1/2$.  When $\xi_1=\cdots=\xi_t$,
it is also true with $n_0=4t+1$ and any choice of $c$ with
$c<1/(4t^2)$ by \cite[Thm.\ A]{RW}.  In general however it is
false and the best that one can hope is an exponent of
approximation of the form $cn^{1/t}$ instead of $cn$ (see
\cite[Prop.\ 10.2]{RW}).  Our main goal here is to show that the
above heuristic statement is true under the restriction that
$\xi_1,\dots,\xi_t$ belong to a field of transcendence degree one
over $\bQ$.  Corresponding values that we find for $n_0$ and $c$
are $n_0=4Dt^2$ and $c=1/(4Dt^3)$ where $D$ denotes the degree of
an algebraic curve of $\bC^t$ defined over $\bQ$ and passing
through the point $(\xi_1,\dots,\xi_t)$.  Although we don't know
if $c$ really requires such dependence in $D$ and $t$, we can
improve its value to $c=1/(4Dt^2)$ when $\xi_1,\dots,\xi_t$ are
distinct. When they further satisfy a recurrence relation of the
form $\xi_{i+1} = a\xi_i+b$ for $i=1,\dots,t-1$, for some $a,b\in
\bQ$ with $a\neq 0,-1$ and $(a,b)\neq (1,0)$, we can even take
$n_0=4t$ and $c=1/(4t)$.

As a bi-product of this work, we obtain a version of Gel'fond's
transcendence criterion expressed in terms of polynomials of bounded
degree taking small values on a fixed sequence of points in
arithmetic progression with rational difference or in geometric
progression with rational ratio.  This new criterion was our
original motivation in writing the present paper, and we hope to
extend its scope in future work.

In the next section, we state our main results in the more general
setting of the previous joint work \cite{RW} with M.~Waldschmidt.
This means that, in order to cover at once the case of approximation
to real, complex or $p$-adic numbers, we replace the field $\bQ$
with a number field $K$, and the field $\bR$ with the completion of
$K$ at some place $w$.  Our strategy for proving these results
follows that of \cite{RW} and is again based on the general method
of Davenport and Schmidt \cite{DS}.  It is briefly described in the
next section.

\section{Main results and notation}
 \label{results}

Throughout this paper, we fix an algebraic extension $K$ of $\bQ$,
a non-trivial place $w$ of $K$, and an algebraic closure $\Kbar$
of $K$. We denote by $d$ the degree $[K:\bQ]$ of $K$ over $\bQ$,
by $\cM$ the set of all non-trivial places of $K$ and by
$\cM_\infty$ the subset of $\cM$ consisting of all Archimedean
places of $K$. For each $v\in\cM$, we denote by $K_v$ the
completion of $K$ at $v$, and by $d_v$ the local degree
$[K_v:\bQ_v]$ where $\bQ_v$ stands for the topological closure of
$\bQ$ in $K_v$.  We also normalize the absolute value $|\ |_v$ of
$K_v$ by asking that, if $v$ is above a prime number $p$ of $\bQ$,
we have $|p|_v=p^{-d_v/d}$ and that, if $v$ is an Archimedean
place, we have $|x|_v=|x|^{d_v/d}$ for any $x\in\bQ$.  With this
convention, the product formula reads $\prod_\vinM |a|_v=1$ for
each $a\in K^\times$, where $K^\times$ stands for the
multiplicative group of $K$.

In order to state our main results of approximation, we also
precise the following notions.  As in the introduction, we define
the {\em height} $\rmH(\alpha)$ of an algebraic number $\alpha$ to
be the largest absolute value of the coefficients of the
irreducible polynomial of $\alpha$ over $\bZ$. The {\em rank} of a
prime ideal $\fp$ of $K[x_1,\dots,x_t]$ is the largest integer
$r\ge 0$ for which there is a strictly increasing chain of $r+1$
prime ideals in $K[x_1,\dots,x_t]$ ending with $\fp$, and its {\em
degree} is defined to be the degree of the corresponding
homogeneous prime ideal of $K[x_0,x_1,\dots,x_t]$. In geometrical
terms, if $V$ denotes the set of zeros of $\fp$ in $\Kbar^t$, then
the rank of $\fp$ is $t-\dim(V)$ and its degree is the degree of
the Zariski closure of $V$ in $\bP^t(\Kbar)$.

\begin{theorem}
 \label{results:thm1}
Let $n$ and $t$ be positive integers, and let $\xi_1,\dots,\xi_t$ be
elements of $K_w$ for which the field $K(\xi_1, \dots, \xi_t)$ has
transcendence degree at most one over $K$.  Let $s$ be the number of
distinct points among $\xi_1,\dots,\xi_t$, let $m$ be the largest
integer for which the sequence $\xi_1,\dots,\xi_t$ contains a point
$\xi$ repeated $m$ times, let $\fp$ be a prime ideal of
$K[x_1,\dots,x_t]$ of rank $t-1$ whose elements all vanish at
$(\xi_1,\dots,\xi_t)$, and let $D$ be the degree of $\fp$. Assume
that
\begin{equation*}
 n \ge 4Dst
 \et
 [K(\xi_i):K] \ge \frac{n}{Dt}
 \quad
 (1\le i\le t).
\end{equation*}
Then, there are infinitely many algebraic numbers $\alpha\in\Kbar$
which, over $K$, have degree $n$ and admit distinct conjugates
$\alpha_1,\dots,\alpha_t$ in $K_w$ satisfying
\begin{equation}
 \label{results:formula1}
 \max_{1\le i\le t}|\xi_i-\alpha_i|_w
 \le \rmH(\alpha)^{-n/(4dDmst)}.
\end{equation}
\end{theorem}

In the case where $\xi_1=\cdots=\xi_t$, we have $s=1$ and $m=t$,
and, as the point $(\xi_1,\dots,\xi_t)$ lies on a rational line,
we can take $D=1$.  Then the above result becomes essentially
\cite[Thm.\ A]{RW}.  In this case, \cite[Prop.\ 10.1]{RW} shows
that the exponent of approximation $n/(4dt^2)$ in
\eqref{results:formula1} is best possible up to the numerical
factor $1/4$.

In the case where $\xi_1,\dots,\xi_t$ are distinct elements of
$K_w$ which all belong to $K+K\xi_1$, we can take $D=m=1$ and
$s=t$, and then the exponent of approximation in
\eqref{results:thm1} again becomes $n/(4dt^2)$.  The following
result provides a special case where the factor $t^2$ in the
denominator can be replaced by $t$.

\begin{theorem}
 \label{results:thm2}
Let $n$ and $t$ be positive integers, and let $\xi_1,\dots,\xi_t$
be elements of $K_w$ which, for some polynomial $\linpol(T)$ of
$K[T]$ of degree one, satisfy the recurrence relation
$\xi_{i+1}=\linpol(\xi_i)$ for $i=1,\dots,t-1$.  Assume moreover
that $\linpol^i(T)\neq T$ for $i=1,\dots,n$, where $\linpol^i$
denotes the $i$-th iterate of $\linpol$, and that
\begin{equation*}
 n \ge 4t
 \et
 [K(\xi_1):K] \ge \frac{n}{t}.
\end{equation*}
Then, there are infinitely many algebraic numbers $\alpha\in\Kbar$
which, over $K$, have degree $n$ and admit distinct conjugates
$\alpha_1, \dots, \alpha_t$ in $K_w$ satisfying
\begin{equation*}
\max_{1\le i\le t}|\xi_i-\alpha_i|_w
  \le \rmH(\alpha)^{-n/(4dt)}.
\end{equation*}
\end{theorem}

Our last main result is the following version of Gel'fond's
transcendence criterion where, for a place $\vinM$ and a polynomial
$Q=\sum_{i=0}^n a_iT^i \in K_v[T]$, we define $\|Q\|_v=\max_{0\le
i\le n} |a_i|_v$.

\begin{theorem}
 \label{results:thm3}
Let $n$ and $t$ be positive integers with $n\ge 4t$, and let
$\xi_1,\dots,\xi_{n+1}$ be elements of $K_w$. Suppose that there
exists a non-zero element $\gamma$ of $K$ such that either we have
$\xi_{i+1}=\gamma+\xi_i$ for $i=1,\dots,n$ (additive case), or we
have $\xi_{i+1}=\gamma\xi_i$ for $i=1,\dots,n$ and $\gamma^i\neq 1$
for $i=1,\dots,2n$ (multiplicative case). Assume moreover that, for
each sufficiently large real number $Y$, there exists a non-zero
polynomial $Q\in K[T]$ of degree at most $n$ which satisfies
$\|Q\|_v\le 1$ for each place $v$ of $K$ distinct from $w$ and also
\begin{equation}
 \label{results:formula2}
 \|Q\|_w \le Y
 \et
 \max_{t+1\le i\le n+1} |Q(\xi_i)|_w  \le Y^{-(4t)/(n+1-4t)}.
\end{equation}
Then, $\xi_1,\dots,\xi_{n+1}$ are algebraic over $K$ of degree
strictly less than $n/t$.
\end{theorem}

Note that the above statement is false if we replace the exponent
$(4t)/(n+1-4t)$ in \eqref{results:formula2} by any exponent
smaller than $t/(n+1-t)$ because Dirichlet box principle shows
that the hypotheses then become satisfied for any choice of
$\xi_1,\dots,\xi_t\in K_w$ (see \cite[\S3]{RW}).

\begin{outline}
For the problem of approximation, instead of looking for just one
polynomial of $K[T]$ of degree at most $n$ taking small values at
the given numbers $\xi_1,\dots, \xi_t$, we follow the idea of
Davenport and Schmidt in \cite{DS} and look for $n+1$ linearly
independent polynomials with this property. Then, it is easy to
build out of them an irreducible polynomial of $K[T]$ of degree $n$
taking small values at $\xi_1, \dots, \xi_t$ while some of its
derivative at the same points are large, so that the new polynomial
has distinct roots $\alpha_1, \dots, \alpha_t$ which are
respectively close to $\xi_1, \dots, \xi_t$, as required. The
precise estimates needed for this are established in \S\ref{app}.

Let $E_n$ denote the $K$-vector space of polynomials of $K[T]$ of
degree at most $n$, and let $g\colon E_n\times E^*_n\to K$ be a
non-degenerate $K$-bilinear form, where $E^*_n$ is any fixed
$K$-vector space with the same dimension $n+1$ as $E_n$.   In order
to produce families of $n+1$ linearly independent polynomials in
$E_n$ as wanted, we use adelic geometry of numbers, asking that the
last minimum of certain adelic convex bodies attached to $E_n$ is at
most one. This is equivalent to asking that slight dilations of
their dual convex bodies, attached to $E^*_n$, have their first
minimum greater than one or, more simply, that they contains no
non-zero element of $E^*_n$.  Precise definitions and relevant
results are given in \S\ref{geom}.

For the approximation results, the choice of the non-degenerate
bilinear form is irrelevant, and we use a standard bilinear form
$\varphi\colon E_n\times K^{n+1} \to K$. Then, there is no useful
interpretation for the elements of the dual convex bodies in
$K^{n+1}$.  However, when the points $\xi_1,\dots,\xi_t$ form an
arithmetic progression with non-zero difference in $K$ or a
geometric progression with non-torsion ratio in $K^\times$, we
construct in \S\ref{bil} special ``translation-invariant'' bilinear
forms $g\colon E_n\times E_n\to K$ for which the dual convex bodies
have the same form as the original ones except that the points
$\xi_1, \dots, \xi_t$ are replaced by the next $n+1-t$ points
$\xi_{t+1},\dots,\xi_{n+1}$ in the corresponding arithmetic or
geometric progression.  Showing that the dual convex bodies contain
no non-zero element of $E_n$ for arbitrarily large values of the
parameters then translates into a version of Gel'fond's criterion in
degree $n$ which is Theorem \ref{results:thm3} above. The reader not
interested in this criterion may skip \S\ref{bil}, while the reader
only interested in it may skip \S\ref{app}. The difficulty of
finding appropriate bilinear form for other choices of points
appears to be an obstacle for extending the criterion to more
general situations.

In \S\ref{mainprop}, we apply the above mentioned principles to
reduce the proof of our main results to the statement that a certain
sequence of adelic convex bodies indexed by a real parameter $X\ge
1$ contains no non-zero element of $K^{n+1}$ for arbitrarily large
values of $X$.  The proof of the latter proceeds by contradiction.
It is done by extending the arguments of \cite[\S\S6--8]{RW} to the
present more general context.  The goal, like in \cite{DS}, is to
replace the sequence of convex bodies by a sequence of polynomials
taking small values at one fixed point and then to derive a
contradiction using an appropriate version of Gel'fond's criterion.
Here we use a version of Gel'fond's criterion for algebraic curves
extending both \cite[Thm.\ 4.2]{RW} and \cite[Thm.\ 2b]{DS}, which
we prove in an appendix.  We apply it to the point
$(\xi_1,\dots,\xi_t)$ and to an algebraic curve containing that
point.  The polynomials that we need are constructed in
\S\ref{constr}, and estimates for their degree and height are
obtained indirectly using auxiliary polynomials in \S\ref{est}. The
proof is completed in \S\ref{proof}.
\end{outline}

\begin{notation}
In the sequel, we use the same notions of heights as in
\cite[\S2]{RW}.

(i) At each place $\vinM$, we define the {\em norm} of an element
$\ux=(x_1,\dots,x_n)$ of $K_v^n$ as its maximum norm
$\|\ux\|_v=\max_{1\le i\le n} |x_i|_v$.  Accordingly, we define
the {\em height} of a point $\ux\in K^n$ by $\rmH(\ux) =
\prod_\vinM \|\ux\|_v$.

(ii) We denote by $E_n$ the vector space over $K$ consisting of
all polynomials $P$ of $K[T]$ with degree at most $n$. We define
the {\em height} $\rmH(P)$ of a polynomial $P\in E_n$ as the
height of its vector of coefficients in $K^{n+1}$. Using the
notation introduced just before the statement of Theorem
\ref{results:thm3}, this is also given by $\rmH(P) = \prod_\vinM
\|P\|_v$.  In view of our notion of height of algebraic numbers,
we also note that, if $\alpha\in\Kbar$ has degree at most $n$ over
$K$ and if $P\in E_n$ is the irreducible polynomial of $\alpha$ in
$K[T]$, then we have $\rmH(\alpha)\le c\rmH(P)^d$ with a constant
$c>0$ depending only on $n$ and $d$ (see \cite[\S2]{RW}).

(iii) Given a place $\vinM$, an integer $m$ with $1\le m\le n$,
and an $m\times n$ matrix $M$ with coefficients in $K_v$, we
define $\|M\|_v$ as the largest absolute value of the minors of
order $m$ of $M$. When $M$ has coefficients in $K$, we define its
{\em height} by $\rmH(M)=\prod_\vinM \|M\|_v$.

(iv) The {\em height} of a subspace $V$ of $E_n$ of dimension
$m\ge 1$ is defined as the height of any $m\times (n+1)$ matrix
$M$ whose rows are the vectors of coefficients of a basis of $V$
over $K$. In particular, if $P$ is a non-zero element of
$E_{n-m+1}$ for some integer $m\ge 1$, and if $V=P\cdot E_{m-1}$
is the subspace of $E_n$ consisting of all products $PQ$ with
$Q\in E_{m-1}$, then Proposition 5.2 of \cite{RW} gives
$c^{-1}\rmH(P)^m\le \rmH(V) \le c\rmH(P)^m$ with a constant $c>0$
which depends only on $n$.

\end{notation}

%
%

\section{Adelic geometry of numbers}
\label{geom}

Let $E$ be a vector space over $K$ of finite dimension $m\ge 1$
(in practice, this will be the space $E_n$ of polynomials of
$K[T]$ of degree at most $n$). For each place $v$ of $K$, we put
$E_v=K_v\otimes_K E$.  We also put $E_\bA=K_\bA\otimes_K E$, where
$K_\bA$ denotes the adele ring of $K$ (see \S14 of \cite{Ca}). We
define on these spaces the natural topology for which any
$K$-linear isomorphism $\psi\colon E\to K^m$ extends by linearity
to a $K_v$-linear homeomorphism $\psi_v\colon E_v\to K_v^m$ for
each $\vinM$, and to a $K_\bA$-linear homeomorphism
$\psi_\bA\colon E_\bA\to K_\bA^m$.  We identify $E$ as a
sub-$K$-vector space of each of these spaces under the natural
embeddings $E\hookrightarrow E_v$ and $E\hookrightarrow E_\bA$
mapping a point $P$ of $E$ to $1\otimes P$.  Then, $E_\bA$ is a
locally compact abelian group and $E$ is a discrete subgroup of
$E_\bA$ with compact quotient $E_\bA/E$.  We equip $E_\bA$ with
the unique Haar measure, denoted $\vol$, for which the quotient
$E_\bA/E$ has measure $1$.

As $K_\bA$ is a topological subring of $\prod_\vinM K_v$, we may
also view $E_\bA$ as a topological subspace of $\prod_\vinM E_v$.
We define a {\it convex body} of $E_\bA$ (or simply of $E$) to be
a compact neighborhood of $0$ in $E_\bA$ of the form
$\cC=\prod_\vinM \cC_v$ where, for each $\vinM$, each
$P,Q\in\cC_v$ and each $a,b\in K_v$, we have $aP+bQ \in \cC_v$
provided that $|a|_v+|b|_v\le 1$ if $v\in\cM_\infty$, or that
$\max\{|a|_v,|b|_v\} \le 1$ if $v\notin\cM_\infty$.  Given a
$K$-linear isomorphism $\psi\colon E \to K^m$, this is equivalent
to asking that $\psi_\bA(\cC)$ is a product $\prod_\vinM \cK_v$
where $\cK_v$ is, in the usual sense, a convex body of $K_v^m$
when $v$ is Archimedean, and a free sub-$\cO_v$-module of $K_v^m$
of rank $m$ otherwise, with $\cK_v=\cO_v^m$ for all but finitely many
ultrametric places $v$, where $\cO_v$ denotes the ring of integers
of $K_v$.

For a convex body $\cC=\prod_\vinM \cC_v$ of $E$ and an idele
$\rho=(\rho_v)_\vinM\in K_\bA^\times$ of $K$, we denote by
$\rho\cC$ the product $\prod_\vinM \rho_v\cC_v$.  This is again a
convex body of $E$.  For a positive real number $\lambda$, we also
define $\lambda\cC$ to be the product $\prod_\vinM \rho_v\cC_v$
where $\rho_v=1$ for each $v\in\cM\setminus\cM_\infty$ and where
$\rho_v=\lambda$ for each $v\in\cM_\infty$ (using the natural
topological embedding of $\bR$ into $K_v$ extending the inclusion
of $\bQ$ into $K$). Finally, for $i=1,\dots,m$, we define the
$i$-th minimum of $\cC$, denoted $\lambda_i(\cC)$, to be the
smallest real number $\lambda>0$ such that $\lambda\cC$ contains
at least $i$ linearly independent elements of $E$ over $K$.

It follows from the above that, if $E'$ is another vector space of
dimension $m$ over $K$ and if $\varphi\colon E \to E'$ is a
$K$-linear isomorphism, then the $K_\bA$-linear map
$\varphi_\bA\colon E_\bA\to E'_\bA$ which extends $\varphi$ maps
any convex body $\cC$ of $E$ to a convex body $\cC'$ of $E'$ with
the same volume and the same successive minima.

In this context, the adelic version of Minkowski's second convex
body theorem proved independently by McFeat \cite{Mc} and by
Bombieri and Vaaler \cite[Thm.\ 3]{BV} reads as follow.

\begin{proposition}
 \label{geom:minkowski}
Let $\cC$ be an adelic convex body of $E$ and let
$\lambda_1,\dots,\lambda_m$ denote its successive minima.  Then,
we have $(\lambda_1\cdots\lambda_m)^d\vol(\cC) \le 2^{md}$.
\end{proposition}

We will also need the following version of Mahler's duality
principle.

\begin{proposition}
 \label{geom:prop:duality}
Let $E^*$ be another vector space over $K$ of dimension $m$, let
$g\colon E\times E^*\to K$ be a non-degenerate $K$-bilinear form,
and let $\cC$ be an adelic convex body of $E$.  For each place $v$
of $K$, define
\begin{equation*}
 \Cg_v = \{ y\in E^*_v \,;\, |g_v(x,y)|_v \le 1
                               \mbox{ for each  $x\in\cC_v$} \},
\end{equation*}
where $g_v\colon E_v\times E^*_v\to K_v$ denotes the
$K_v$-bilinear form which extends $g$.  Then, $\Cg=\prod_\vinM
\Cg_v$ is an adelic convex body of $E^*$.  Moreover, if
$\lambda_1,\dots,\lambda_m$ denote the successive minima of $\cC$
and $\lambdag_1,\dots,\lambdag_m$ those of $\Cg$, then we have $1
\le \lambda_i\lambdag_{m+1-i} \le c_1$ for $i=1,\dots,m$ with a
constant $c_1\ge 1$ depending only on $K$ and $m$.
\end{proposition}

We refer to $\Cg$ as the {\it dual} of $\cC$ with respect to $g$.

\begin{proof}
This follows from Lemma 3.1 (ii) and  Theorem 3.7 of \cite{Bu}, in
the case where $E=E^*=K^m$ and where $g$ is the usual bilinear
form $\theta\colon K^m\times K^m\to K$ given by $\theta(\ux,\uy) =
\sum_{i=1}^m x_iy_i$ for each $\ux=(x_1,\dots,x_m)$ and each
$\uy=(y_1,\dots,y_m)$ in $K^m$ (see also Theorem 1 of \cite{Du}
for the case $K=\bQ$). To deduce the general case, choose
$K$-linear isomorphisms $\psi\colon E \to K^m$ and $\psi^*\colon
E^*\to K^m$ satisfying $g(P,Q)=\theta(\psi(P),\psi^*(Q))$ for each
$(P,Q)\in E\times E^*$.  Put $\cK=\prod_\vinM \psi_v(\cC_v)$ and
$\Ktheta=\prod_\vinM \psi_v^*(\cC^*_v)$. Then, $\Ktheta$ is the
dual of $\cK$ with respect to $\theta$ and so, by Lemma 3.1 of
\cite{Bu}, $\Ktheta$ is a convex body of $K^m$ which in turn
implies that $\Cg$ is a convex body of $E^*$. Moreover, the
successive minima of $\cK$ and $\Ktheta$ being respectively the
same as those of $\cC$ and $\Cg$, Theorem 3.7 of \cite{Bu} gives
$1 \le \lambda_i\lambdag_{m+1-i} \le c_1$ for $i=1,\dots,m$ with
an explicit constant $c_1=c_1(K,m)$.
\end{proof}

\begin{remark}
With the notation of Proposition \ref{geom:prop:duality}, if a
point $P\in E$ belongs to the interior of $\cC$ and if $Q\in E^*$
belongs to $\Cg$, then the element $g(P,Q)$ of $K$ satisfies
$|g(P,Q)|_v < 1$ for each $v\in\cM_\infty$ and $|g(P,Q)|_v \le 1$
for each $v\in\cM\setminus\cM_\infty$.  This gives $\prod_\vinM
|g(P,Q)|_v <1$ and so $g(P,Q)=0$ by virtue of the product formula.
\end{remark}

We also recall the statement of the strong approximation theorem
(see \cite[Thm.\ 3, p.\ 440]{Ma} or \cite[\S15]{Ca}).

\begin{lemma}
 \label{geom:lemma1}
There exists a constant $c_2>0$ depending only on $K$ with the
following property.  Let $\cS$ be a finite set of places of $K$
and, for each $v\in \cS$, let $\theta_v$ be an element of $K_v$
and let $\epsilon_v$ be a positive real number.  Assume that
$\prod_{v\in\cS} \epsilon_v \ge c_2$.  Then, there exists an
element $a$ of $K$ satisfying $|a-\theta_v|_v \le \epsilon_v$ for
each $v\in \cS$, and $|a|_v\le 1$ for each $v\in\cM\setminus\cS$.
\end{lemma}

As a first application, we note the following simple consequence.

\begin{proposition}
 \label{geom:prop:renormalization}
Let $\cC=\prod_\vinM \cC_v$ be an adelic convex body of $E$, and
let $\rho=(\rho_v)_\vinM \in K_\bA^\times$ be an idele of $K$.
Define the {\em content} of $\rho$ to be $c(\rho)=\prod_\vinM
|\rho_v|_v$.  Then, for $i=1,\dots,m$, we have
$c_2^{-1}\lambda_i(\cC) \le c(\rho)\lambda_i(\rho\cC) \le c_2
\lambda_i(\cC)$.
\end{proposition}

\begin{proof}
Fix an index $i$ with $1\le i\le m$ and put
$\lambda=\lambda_i(\cC)$, so that $\lambda\cC$ contains $i$
linearly independent elements $P_1,\dots,P_i$ of $E$.  Choose also
a real number $c$ with $c>c_2c(\rho)^{-1}$, and an ultrametric
place $u$ of $K$ with $|\rho_u|_u=1$.  Then, by Lemma
\ref{geom:lemma1}, there exists an element $a$ of $K$ satisfying
$|a|_v\le c^{d_v/d}|\rho_v|_v$ for each $v\in\cM_\infty$,
$|a|_v\le |\rho_v|_v$ for each $v\notin \cM_\infty\cup\{u\}$, and
also $|a-1|_u<1$.  The last condition $|a-1|_u<1$ implies that $a$
is non-zero and gives $|a|_u\le |\rho_u|_u$. Then
$aP_1,\dots,aP_i$ are linearly independent elements of $E$ which
belong to $a\lambda\cC_v \subseteq c\lambda\rho_v\cC_v$ for each
$v\in\cM_\infty$ and belong to $a\cC_v \subseteq \rho_v\cC_v$ for
each $v\notin\cM_\infty$, showing that $\lambda_i(\rho\cC) \le
c\lambda$.  By virtue of the choice of $c$, this proves that
$c(\rho)\lambda_i(\rho\cC) \le c_2 \lambda_i(\cC)$.  The lower
bound for $c(\rho)\lambda_i(\rho\cC)$ follows from this inequality
by applying it to the pair $\rho\cC$ and $\rho^{-1}$ instead of
$\cC$ and $\rho$, upon noting that $c(\rho^{-1})=c(\rho)^{-1}$.
\end{proof}

The last proposition formalizes the construction of Davenport and
Schmidt in \cite[\S 2]{DS}.

\begin{proposition}
 \label{geom:prop:appsim}
Let $\cC=\prod_\vinM \cC_v$ be an adelic convex body of $E$, and
let $\cS$ be a finite set of places of $K$.  For each $v\in\cS$,
choose $P_v\in E_v$ and $\rho_v\in K_v$ satisfying
$\prod_{v\in\cS} |\rho_v|_v \ge c_2 m\lambda_m(\cC)$. Then, there
exists an element $P$ of $E$, satisfying $P-P_v \in \rho_v\cC_v$
for each $v\in \cS$, and $P\in \cC_v$ for each $v\in
\cM\setminus\cS$.
\end{proposition}

\begin{proof}
Defining $P_v=0$ and $\rho_v=1$ for each Archimedean place $v$ of
$K$ not in $\cS$, we may assume, without loss of generality, that
$\cM_\infty \subseteq \cS$. Put $\lambda=\lambda_m(\cC)$. By
definition, the convex body $\lambda\cC$ contains a basis
$\{P_1,\dots,P_m\}$ of $E$ over $K$. For each place $v\in\cS$,
this basis is also a basis of $E_v$ over $K_v$ and so we can write
\begin{equation*}
 P_v=\theta_{1,v}P_1+\cdots+\theta_{m,v}P_m
\end{equation*}
with $\theta_{1,v},\dots\theta_{m,v}\in K_v$.  Define $\epsilon_v
= |\rho_v|_v$ for each $v\in\cS\setminus\cM_\infty$ and
$\epsilon_v = (m\lambda)^{-d_v/d}|\rho_v|_v$ for each
$v\in\cM_\infty$. Since we have $\prod_{v\in\cS} \epsilon_v \ge
c_2$, Lemma \ref{geom:lemma1} provides, for $i=1,\dots,m$, an
element $a_i$ of $K$ satisfying $|a_i|_v\le 1$ for any
$v\in\cM\setminus\cS$, and $|a_i-\theta_{i,v}|_v\le \epsilon_v$
for any $v\in \cS$. We claim that the polynomial
$P=a_1P_1+\cdots+a_mP_m$ has all the required properties.  First
of all, for each $v\in\cM_\infty$ and for $i=1,\dots,m$, we have
$|a_i-\theta_{i,v}|_v \le |(m\lambda)^{-1}\rho_v|_v$ and
$P_i\in\lambda\cC_v$, so that all products $(a_i-\theta_{i,v})P_i$
belong to $m^{-1}\rho_v\cC_v$ and their sum $P-P_v$ belongs to
$\rho_v\cC_v$.  For each $v\in\cS\setminus\cM_\infty$, we have
instead $|a_i-\theta_{i,v}|_v \le |\rho_v|_v$ and $P_i\in\cC_v$
for each $i$, so that $\rho_v\cC_v$ contains all products
$(a_i-\theta_{i,v})P_i$ and also their sum $P-P_v$. Finally, for
each of the remaining ultrametric places $v\in \cM\setminus\cS$,
we have $a_iP_i\in a_i\cC_v \subseteq \cC_v$ for each $i$, and so
$P\in\cC_v$.
\end{proof}

%
%

\section{Approximation}
 \label{app}

Let $n$ and $t$ be integers with $1\le t\le n$, and let $(\xi_1,
\dots, \xi_t)$ be a point of $K_w^t$. We denote by
$\eta_1,\dots,\eta_s$ the distinct elements of the sequence
$(\xi_1,\dots,\xi_t)$ and, for each $i=1,\dots,s$, we denote by
$m_i$ the number of times that $\eta_i$ appears in this sequence.
The constants $c_3,c_4,\dots$ that appear below, as well as the
implied constants in the symbols $\ll$ and $\gg$ depend only on
$K$, $w$ and the above quantities.

As mentioned in \S\ref{results}, we denote by $E_n$ the vector
space over $K$ consisting of all polynomials of $K[T]$ of degree
$\le n$. Then, for each $\vinM$, the $K_v$-vector space $E_{n,v} =
K_v\otimes_k E_n$ identifies itself with the space of polynomials
of $K_v[T]$ of degree $\le n$.  For each pair of real numbers
$X,Y\ge 1$, we define $\cC(X,Y)$ to be the convex body of $E_n$
whose component $\cC_v(X,Y)$ at each place $\vinM$ distinct from
$w$ consists of all polynomials $P$ of $E_{n,v}$ with $\|P\|_v\le
1$, and whose component $\cC_w(X,Y)$ at $w$ consists of all
polynomials $P$ of $E_{n,w}$ with
\begin{equation}
 \label{constr:convex}
 \|P\|_w\le X
 \et
 \max_{1\le i\le s}
 \Big(
   \max_{0\le j<m_i} |P^{(j)}(\eta_i)|_w
 \Big)
 \le Y^{-1},
\end{equation}
where $P^{(j)}$ stands for the $j$-th derivative of $P$.  Our goal
in this section is to prove the following result of approximation
which in a sense extends Lemma 9.1 of \cite{RW}.

\begin{proposition}
 \label{app:prop1}
Let $X,Y$ be real numbers with $X,Y\ge 1$ and let $\lambda$ denote
the $(n+1)$-th minimum of the convex body $\cC(X,Y)$ of $E_n$.
Then there exists an irreducible polynomial $P\in K[T]$ of degree
$n$ and height at most $c_3\lambda X$ which, for $i=1,\dots,s$,
admits at least $m_i$ roots in the closed disk of $K_w$ of radius
$c_4(XY)^{-1/m_i}$ centered at $\eta_i$, without vanishing at
$\eta_i$.
\end{proposition}

\begin{proof}
Choose a finite place $u$ of $K$ with $u\neq w$, and a
uniformizing parameter $\pi\in\cO_u$ for $u$.  Define $P_u=T+\pi$
and $\rho_u=\pi^2$.  By Proposition \ref{geom:prop:appsim}, there
is a constant $c_5>0$ such that, for any choice of $P_w\in
E_{n,w}$ and $\rho_w\in K_w$ with $|\rho_w|_w\ge c_5 \lambda$,
there is a polynomial $P\in E_n$ which satisfies $P\in\cC_v(X,Y)$
for each $v \in \cM\setminus\{u,w\}$, and $P-P_v\in
\rho_v\cC_v(X,Y)$ for each $v \in \{u,w\}$. The condition at $u$
reads
\begin{equation*}
 \| P - (T^n+\pi) \|_u \le |\pi|_u^2.
\end{equation*}
By virtue of Eisenstein's criterion, it implies that such a
polynomial $P$ is irreducible over $K_u$ and so is irreducible
over $K$. It also gives $\|P\|_u\le 1$.  Since, for each $v \notin
\{u,w\}$, the condition $P\in\cC_v(X,Y)$ means $\|P\|_v\le 1$, we
deduce that
\begin{equation}
 \label{app:height}
 \rmH(P) \le \|P\|_w.
\end{equation}
Choose $\rho_w$ to be an element of $K_w$ of smallest norm with
$|\rho_w|_w\ge c_5\lambda$.  Then, we have $|\rho_w|_w\ll
\lambda$, and the last condition $P-P_w\in \rho_w\cC_w(X,Y)$ leads
to
\begin{equation}
 \label{app:condw}
  \begin{aligned}
  \|P-P_w\|_w &\ll \lambda X, \\
  |(P-P_w)^{(j)}(\eta_i)|_w &\ll \lambda Y^{-1}
  \quad
  (1\le i\le s,\ 0\le j< m_i).
  \end{aligned}
\end{equation}
 We look for a polynomial $P_w$ of the form
\begin{equation*}
 P_w(T) = a \prod_{i=1}^s \prod_{j=1}^{m_i} (T-\eta_i-jz_i)
\end{equation*}
with $a,z_1,\dots,z_s \in K_w\setminus\{0\}$. To choose the latter
parameters, we note that there exists a constant $c_6$ with $0<
c_6< (n!)^{-1}$ such that any polynomial $Q\in E_{n,w}$ which
satisfies
\begin{equation*}
 \label{app:delta}
 \Big\| Q(T) - \prod_{j=1}^\ell (T-j) \Big\|_w \le c_6
\end{equation*}
for some integer $\ell$ with $1\le \ell \le n$ admits at least
$\ell$ distinct roots of norm at most $|\ell|_w+1$ in $K_w$.
Putting
\begin{equation*}
 \sigma_i
 = a z_i^{m_i} \prod_{k\neq i} (\eta_i-\eta_k)^{m_k}
 \quad
 (1\le i \le s),
\end{equation*}
where the product extends to all integers $k=1,\dots,s$ with
$k\neq i$, we find
\begin{equation*}
 \sigma_i^{-1} P_w(z_iT+\eta_i) - \prod_{j=1}^{m_i} (T-j)
 =
 \Big( \prod_{k\neq i} \prod_{j=1}^{m_k}
         \Big(
         1 + \frac{z_iT-jz_k}{\eta_i-\eta_k}
         \Big) -1
 \Big)
 \prod_{j=1}^{m_i} (T-j)
\end{equation*}
so that
\begin{equation}
 \label{app:diff}
 \Big\|
 \sigma_i^{-1} P_w(z_iT+\eta_i) - \prod_{j=1}^{m_i} (T-j)
 \Big\|_w
 \le
 c_7 \max_{1\le k\le s} |z_k|_w
\end{equation}
for some constant $c_7\ge 1$.  We now fix $z_1,\dots,z_s\in K_w$
of maximal absolute value with
\begin{equation*}
 |z_i|_w \le \frac{c_6}{2c_7}(XY)^{-1/m_i}  \quad (1\le i\le s).
\end{equation*}
Then the right hand side of \eqref{app:diff} is bounded above by
$c_6/2$ and using the hypotheses \eqref{app:condw} on $P$ we find,
for each $i=1,\dots,s$,
\begin{equation}
 \label{app:diffP}
 \begin{aligned}
 \Big\|
 \sigma_i^{-1} P(z_iT+\eta_i) - \prod_{j=1}^{m_i} (T-j)
 \Big\|_w
 &\le
 \frac{c_6}{2}
   + |\sigma_i|_w^{-1}\|P(z_iT+\eta_i)-P_w(z_iT+\eta_i)\|_w \\
 &=
 \frac{c_6}{2}
   + |\sigma_i|_w^{-1}
     \max_{0\le j\le n}
     \Big| \frac{1}{j!} (P-P_w)^{(j)}(\eta_i) \Big|_w |z_i|_w^j\\
 &\le
 \frac{c_6}{2} +c_8\lambda |a|_w^{-1} X,
 \end{aligned}
\end{equation}
with a constant $c_8>0$.  Finally, we fix $a\in K_w$ of minimal
absolute value with
\begin{equation*}
 |a|_w \ge \frac{2c_8}{c_6} \lambda X.
\end{equation*}
Then the left hand side of \eqref{app:diffP} is at most $c_6$ and
accordingly the polynomial $P$ admits at least $m_i$ distinct
roots in the ball of $K_w$ of radius $(|m_i|_w+1)|z_i|_w$ centered
at $\eta_i$. Moreover, substituting $T=0$ in \eqref{app:diffP}
provides $|\sigma_i^{-1} P(\eta_i) \pm m_i! |_w  \le c_6 <
|m_i!|_w$, and so we must have $P(\eta_i)\neq 0$. The choice of
$a$ also gives $\|P_w\|_w\ll \lambda X$. Combining this with
\eqref{app:height} and \eqref{app:condw}, we deduce that
$\rmH(P)\ll \lambda X$. Thus $P$ has all the required properties.
\end{proof}

%
%

\section{Invariant bilinear forms}
\label{bil}

Assume that the points $\xi_1,\dots,\xi_t$ introduced in the
previous section \S\ref{app} come from a sequence
$\xi_1,\dots,\xi_{n+1}$ of $n+1$ distinct elements of $K_w$ which
is either an arithmetic progression with difference $\gamma\in
K^\times$ (the additive case), or a geometric progression with
ratio $\gamma\in K^\times$ satisfying $\xi_1\neq 0$ and
$\gamma^i\neq 1$ for $i=1,\dots,2n$ (the multiplicative case).

Then, for each pair of real numbers $X,Y\ge 1$, the convex body
$\cC(X,Y)$ of $E_n$ introduced in the previous section is the
product $\prod_\vinM \cC_v(X,Y)$ where
\begin{equation}
 \label{bil:C}
 \cC_w(X,Y)
 =
 \{ P\in E_{n,w} \,;\, \|P\|_w \le X \mbox{ and }
        \max_{1\le i\le t} |P(\xi_i)|_w \le Y^{-1}\},
\end{equation}
and where, for $v\neq w$, the component $\cC_v(X,Y)$ consists of
all polynomials $P$ of $E_{n,v}$ with $\|P\|_v\le 1$.  Similarly,
we define another convex body $\Cbar(X,Y) = \prod_\vinM
\Cbar_v(X,Y)$ by putting
\begin{equation}
 \label{bil:Cbar}
 \Cbar_w(X,Y)
 =
 \{ Q\in E_{n,w} \,;\, \|Q\|_w \le Y \mbox{ and }
        \max_{t+1\le i\le n+1} |P(\xi_i)|_w \le X^{-1}\},
\end{equation}
and $\Cbar_v(X,Y)=\cC_v(X,Y)$ for every $v\neq w$.  Our goal is to
show that these convex bodies $\cC(X,Y)$ and $\Cbar(X,Y)$ are
essentially dual to each other with respect to the bilinear form
$g$ constructed by the following proposition.

\begin{proposition}
 \label{bil:prop1}
Let $\gamma\in K^\times$ be as above (with $\gamma^i\neq 1$ for
$i=1,\dots,2n$ in the multiplicative case).  For each integer
$i\ge 0$, we define $\gamma_i=i\gamma$ in the additive case, and
$\gamma_i=\gamma^i$ in the multiplicative case.  For each $x\in
K$, we also denote by $\tau_x\colon E_n \to E_n$ the linear map
given by $\tau_x(P(T)) = P(x+T)$ in the additive case, and by
$\tau_x(P(T)) = P(xT)$ in the multiplicative case. Then, in each
case, there exist elements $g_{ij}$ of $K$ for $0\le i \le j\le n$
such that the bilinear form $g\colon E_n\times E_n \to K$ given by
\begin{equation}
 \label{bil:form}
   g(P,Q) = \sum_{0\le i\le j\le n} g_{ij} P(\gamma_i)Q(\gamma_j)
\end{equation}
is non-degenerate and satisfies, for any $x\in K$ and any $P,Q\in
E_n$,
\begin{equation}
 \label{bil:invariance}
  g(\tau_xP,\tau_xQ) =
  \begin{cases}
    g(P,Q)     &\text{in the additive case,} \\
    x^n g(P,Q) &\text{in the multiplicative case.}
  \end{cases}
\end{equation}
\end{proposition}

\begin{proof}
The hypothesis on $\gamma$ ensures that the points $\gamma_0,
\gamma_1, \dots, \gamma_{2n}$ are all distinct.  In particular,
the first $n+1$ of them are distinct and so there exists a unique
choice of elements $a_0,a_1,\dots,a_{n+1}$ of $K$ with $a_{n+1}=1$
such that
 \begin{equation}
   \label{bil:relation}
   \sum_{i=0}^{n+1} a_i P(\gamma_i) = 0
 \end{equation}
for any $P\in E_n$.

Fix temporarily $P\in E_n$.  In the additive case, we put $\rho=1$
and $\tilde{P}(T) = P((n+1)\gamma-T)$.  In the multiplicative
case, we put $\rho=\gamma^n$ and $\tilde{P}(T) = T^n
P(\gamma^{n+1}T^{-1})$. Then $\tilde{P}$ belongs to $E_n$ and the
formula \eqref{bil:relation} applied to $\tilde{P}$ becomes
 \[
 \sum_{i=0}^{n+1} \rho^i a_i P(\gamma_{n+1-i}) = 0.
 \]
 From this we deduce that
 \begin{equation}
   \label{bil:relationbis}
   \sum_{i=0}^{n+1} \rho^{n+1-i}a_{n+1-i} P(\gamma_i) = 0
 \end{equation}
for any $P\in E_n$ and so, although we will not need it, we get
$\rho^{n+1-i}a_{n+1-i} = a_0a_i$ for $i=0,1,\dots,n+1$.

Let $g\colon E_n\times E_n\to K$ be the bilinear form given by
\eqref{bil:form} for the choice of coefficients
\begin{equation}
 \label{bil:coeff}
   g_{ij} = \rho^{n-j} a_{n+1+i-j} \quad (0\le i\le j\le n).
\end{equation}
Since $g_{ii}\neq 0$ for $i=0,1,\dots,n$, this bilinear form is
non-degenerate.  Moreover, for any pair of polynomials $P,Q\in
E_n$, we find, using \eqref{bil:relation}, \eqref{bil:relationbis}
and \eqref{bil:coeff},
\begin{align*}
 g(\tau_\gamma P, \tau_\gamma Q)
 &= \sum_{0\le i\le j\le n-1} g_{ij} P(\gamma_{i+1})Q(\gamma_{j+1})
    + \sum_{0\le i\le n} a_{i+1}P(\gamma_{i+1})
      Q(\gamma_{n+1}) \\
 &= \sum_{1\le i\le j\le n} \rho g_{ij} P(\gamma_i)Q(\gamma_j)
    - a_0 P(\gamma_0) Q(\gamma_{n+1}) \\
 &= \rho g(P,Q)
    - P(\gamma_0) \sum_{0\le j\le n+1} \rho^{n+1-j} a_{n+1-j}
                        Q(\gamma_j) \\
 &=\rho g(P,Q).
\end{align*}
By recurrence, we deduce that the formula \eqref{bil:invariance}
holds for $x=\gamma_i$ with $i=0,1,\dots,2n$.  Since these $2n+1$
numbers are distinct, and since $g(\tau_xP,\tau_xQ)$ is, for fixed
$P$ and $Q$, a polynomial in $x$ of degree at most $2n$, the
formula must therefore hold for any $x\in K$.
\end{proof}

\begin{remarks}
In the additive (resp.\ multiplicative) case, the property
\eqref{bil:invariance} expresses an invariance of the bilinear
form $g$ under the additive (resp.\ multiplicative) group of $K$.
One may wonder if similar invariant ``triangular'' forms can be
defined for elliptic curves defined over $K$.  In the present
context, it is interesting to note that the bilinear form $g\colon
E_n\times E_n \to K$ defined in \cite[Lemma 3.3]{RW} in terms of
the derivatives of the polynomials at $0$ possesses both
invariance properties stated in \eqref{bil:invariance}.
\end{remarks}

\begin{remarks}
In the notation of the proof, one finds for $i=0,1,\dots,n$
that $a_i=-P_i(\gamma_{n+1})$ where $P_i$ denotes the element of
$E_n$ which takes the value $1$ at $\gamma_i$ and vanishes at all
other points $\gamma_j$ with $0\le j\le n$ and $j\neq i$, and
therefore
 \[
   a_i = -\prod_{j\neq i}
          \frac{\gamma_{n+1}-\gamma_j}{\gamma_i-\gamma_j}
 \]
where the product extends over all indices $j$ with $0\le j\le n$
with $j\neq i$.  In particular, in the additive case, we find that
$a_i=(-1)^{n+1-i} \binom{n+1}{i}$ is independent of $\gamma$.  We
will not need these explicit formulas here.
\end{remarks}

\medskip
Coming back to the adelic convex bodies $\cC(X,Y)$ and
$\Cbar(X,Y)$, we can now state the result which was alluded to at
the beginning of the section.

\begin{proposition}
 \label{bil:prop2}
There exist ideles $\alpha,\beta \in K_\bA^\times$ such that, for
any choice of real numbers $X,Y\ge 1$, we have
\begin{equation}
 \label{bil:comp}
 \alpha\Cbar(X,Y)\ \subseteq \Cg(X,Y) \subseteq\ \beta\Cbar(X,Y),
\end{equation}
where $\Cg(X,Y)$ denotes the dual of $\cC(X,Y)$ with respect to
the bilinear form $g$ given by Proposition \ref{bil:prop1}.
\end{proposition}

\begin{proof}
Let $g_w\colon E_{n,w}\times E_{n,w} \to K_w$ denote the
$K_w$-bilinear form which extends $g$. The invariance property
\eqref{bil:invariance} of $g$ extends by continuity to $g_w$ for
each $x\in K_w$, with the map $\tau_x\colon E_{n,w}\to E_{n,w}$
defined by the same formula as in Proposition \ref{bil:prop1}.
Taking $x=\xi_1$, this gives
\begin{equation}
 \label{bil:formbis}
 g_w(P,Q)
 =
 \sum_{0\le i\le j\le n} \rho g_{ij} P(\xi_{i+1})Q(\xi_{j+1})
\end{equation}
for any $P,Q\in E_{n,w}$, with $\rho=1$ in the additive case and
$\rho=\xi_1^{-n}$ in the multiplicative case.  From this and the
definitions \eqref{bil:C} and \eqref{bil:Cbar}, we deduce that
$|g_w(P,Q)|_w \le c$ for each $P\in\cC_w(X,Y)$ and each
$Q\in\Cbar_w(X,Y)$, with a constant $c>0$ that is independent of
$X$ and $Y$.  Choosing $\alpha_w\in K_w^\times$ with $|\alpha_w|_w
\le 1/c$ then gives $\alpha_w\Cbar_w(X,Y) \subseteq \Cg_w(X,Y)$.
Conversely, since $\xi_1,\dots,\xi_{n+1}$ are distinct, there
exist, for each $i=1,\dots,n+1$, a unique polynomial $P_i\in
E_{n,w}$ which vanishes at each point $\xi_j$ with $j\neq i$ and
satisfies $P_i(\xi_i)=Y^{-1}$ if $i\le t$, and $P_i(\xi_i)=X$ if
$i> t$. The polynomials $P_1, \dots, P_{n+1}$ so constructed
belong to $\theta\cC_w(X,Y)$ for some constant $\theta\in
K_w^\times$. Then, any $Q\in\Cg_w(X,Y)$ satisfies
$|g_w(\theta^{-1}P_i,Q)|_w\le 1$ for $i=1,\dots,n+1$ which, in
view of \eqref{bil:formbis}, translates into
\begin{equation*}
 \left| \sum_{j=i}^n \rho g_{ij} \theta^{-1} Q(\xi_{j+1}) \right|_w
 \le
 \begin{cases}
  Y      &\mbox{for $i=0,\dots,t-1$,}\\
  X^{-1} &\mbox{for $i=t,\dots,n$.}
 \end{cases}
\end{equation*}
As $g_{ii}\neq 0$ for $i=0,\dots,n$, this implies that
$\Cg_w(X,Y)$ is contained in $\beta_w\Cbar_w(X,Y)$ for a constant
$\beta_w\in K_w^\times$.  For the remaining places $v\neq w$ of
$K$, the components of $\cC(X,Y)$ and $\Cbar(X,Y)$ at $v$ are
independent of $X$ and $Y$, and thus we also have
\begin{equation*}
 \alpha_v\Cbar_v(X,Y)
 \subseteq \Cg_v(X,Y)
 \subseteq \beta_v \Cbar_v(X,Y)
\end{equation*}
for some $\alpha_v, \beta_v\in K_v^\times$ which are independent
of $X$ and $Y$ and can be taken to be $1$ for all but finitely
many places $v$. Then the ideles $\alpha=(\alpha_v)_\vinM$ and
$\beta=(\beta_v)_\vinM$ have the property \eqref{bil:comp}.
\end{proof}

Combining this result with Propositions \ref{geom:prop:duality}
and \ref{geom:prop:renormalization}, we deduce the following.

\begin{corollary}
 \label{bil:cor}
Let $\lambda_j(X,Y)$ and $\lambdabar_j(X,Y)$ denote respectively
the $j$-th minima of $\cC(X,Y)$ and $\Cbar(X,Y)$.  Then, the
products $\lambda_j(X,Y)\lambdabar_{n+2-j}(X,Y)$ with $1\le j\le
n+1$ are bounded above and below by positive constants which are
independent of the choice of $X,Y\ge 1$.
\end{corollary}

%
%

\section{The main Proposition and deduction of the Theorems}
 \label{mainprop}

Let the notation be as in \S\ref{app}, and let $\varphi\colon
E_n\times K^{n+1} \to K$ be the non-degenerate $K$-bilinear form
given by
\begin{equation*}
 \varphi(a_0+a_1T+\cdots+a_nT^n,(y_0,y_1,\dots,y_n))
 =
 a_0y_0+a_1y_1+\cdots+a_ny_n.
\end{equation*}
In \S\ref{app}, we defined a convex body $\cC(X,Y)$ for each pair
of real numbers $X,Y\ge 1$.  Accordingly, we denote by
$\Cphi(X,Y)$ the convex body of $K^{n+1}$ which is dual to
$\cC(X,Y)$ with respect to $\varphi$. For $i=1,\dots,n+1$, we also
denote by $\lambda_i(X,Y)$ and $\lambdaphi_i(X,Y)$ the respective
$i$-th minimum of $\cC(X,Y)$ in $E_n$ and of $\Cphi(X,Y)$ in
$K^{n+1}$. We show in this section how the theorems stated in
\S\ref{results} can be derived from the following proposition
whose proof is postponed to the last section \S\ref{proof}.

\begin{proposition}
 \label{mainprop:prop1}
Assume that we are either in the situation of Theorem
\ref{results:thm1}, in which case we define $\nu=4Dst$, or in the
situation of Theorem \ref{results:thm2}, in which case we define
$\nu=4t$.  Then, there are arbitrarily large values of $X$ such
that $\lambdaphi_1(X,X^{(n+2-\nu)/\nu})>1$.
\end{proposition}

\subsection{Proof of Theorems \ref{results:thm1} and \ref{results:thm2}}
Assume that we are in the situation of Theorem \ref{results:thm1}
or Theorem \ref{results:thm2} and define $\nu$ accordingly as in
Proposition \ref{mainprop:prop1}.  For each pair of real numbers
$X, Y\ge 1$ satisfying $\lambdaphi_1(X,Y)>1$, Proposition
\ref{geom:prop:duality} gives $\lambda_{n+1}(X,Y) \ll
\lambdaphi_1(X,Y)^{-1} \ll 1$ with implied constants which are
independent of $X$ and $Y$, and then Proposition \ref{app:prop1}
shows the existence of an irreducible polynomial $P\in K[T]$ of
degree $n$ and height $\ll X$ which, for $i=1,\dots,s$, admits at
least $m_i$ roots in a closed disk of $K_w$ of radius $\ll
(XY)^{-1/m_i}$ centered at $\eta_i$, without vanishing at
$\eta_i$. If the product $XY$ is sufficiently large, these disks
are disjoint and we deduce that $P$ admits $t$ distinct roots
$\alpha_1,\dots,\alpha_t$ satisfying $0<|\xi_i-\alpha_i|_w \ll
(XY)^{-1/m}$, where $m=\max\{m_1,\dots,m_s\}$.  Moreover, if
$\alpha$ is a root of $P$ in $\Kbar$, then we have $\rmH(\alpha)
\ll \rmH(P)^d \ll X^d$.  Assume from now on that
$Y=X^{(n+2-\nu)/\nu}$. Then the hypothesis gives
$\lambdaphi_1(X,Y)>1$ for arbitrary large values of $X$ and, for
each such $X$, the above provides an algebraic number
$\alpha=\alpha_X \in \Kbar$ which, over $K$, has degree $n$ and
admits distinct conjugates $\alpha_1,\dots,\alpha_t \in K_w$
satisfying $0<|\xi_i-\alpha_i|_w \ll X^{-(n+2)/(m\nu)} \ll
\rmH(\alpha)^{-(n+2)/(dm\nu)}$ for $i=1,\dots,t$.  The conclusion
follows as, by varying $X$, we get infinitely many algebraic
numbers $\alpha$.

\subsection{Proof of Theorem \ref{results:thm3}}
Let the notation and hypotheses be as in Theorem
\ref{results:thm3}. For each pair of real numbers $X,Y\ge 1$,
define $\Cbar(X,Y)$ as in \S\ref{bil} and, for $i=1,\dots,n+1$,
denote by $\lambdabar_i(X,Y)$ the $i$-th minimum of $\Cbar(X,Y)$.
According to Proposition \ref{geom:prop:duality} and Corollary
\ref{bil:cor}, the products
$\lambdaphi_i(X,Y)\lambda_{n+2-i}(X,Y)$ and
$\lambdabar_i(X,Y)\lambda_{n+2-i}(X,Y)$ are bounded below and
above by positive constants which are independent of $X$ and $Y$.
The same is therefore true of the ratios
$\lambdaphi_i(X,Y)/\lambdabar_i(X,Y)$.  In particular there exists
a constant $c>0$ such that $\lambdaphi_1(X,Y)\le
c\lambdabar_1(X,Y)$. Moreover, if $\rho=(\rho_v)_\vinM\in
K_\bA^\times$ is an idele of $K$ satisfying $\rho_v=1$ for each
place $v\neq w$, then, putting $r=|\rho_w|_w$, we find
$\rho\Cbar(X,Y)=\Cbar(r^{-1}X,rY)$ and accordingly Proposition
\ref{geom:prop:renormalization} gives $\lambdabar_i(r^{-1}X,rY)\le
c_2r^{-1}\lambdabar_i(X,Y)$ for $i=1,\dots,n+1$.  In particular,
for a suitable choice of $r\ge 1$, we have
\begin{equation}
 \label{mainprop:lambda}
 \lambdaphi_1(r^{-1}X,rY)
 \le c \lambdabar_1(r^{-1}X,rY)
 \le \lambdabar_1(X,Y),
\end{equation}
independently of $X\ge r$ and $Y\ge 1$.

The hypotheses of Theorem \ref{results:thm3} imply that, for each
choice of $X$ and $Y$ with $1\le X\le Y^{4t/(n+1-4t)}$ and $Y$
sufficiently large, the convex body $\Cbar(X,Y)$ contains a
non-zero element of $E_n$, and so we have $\lambdabar_1(X,Y)\le
1$.  By \eqref{mainprop:lambda}, this gives
$\lambdaphi_1(r^{-1}X,rY) \le 1$ and thus $\Cphi(r^{-1}X,rY)$
contains a non-zero element of $K^{n+1}$ for such choices of $X$
and $Y$. In particular, we deduce that
$\Cphi(X,X^{(n+2-4t)/(4t)})$ contains a non-zero element of
$K^{n+1}$ for each sufficiently large value of $X$.  By
Proposition \ref{mainprop:prop1}, this means that the given points
$\xi_1,\dots,\xi_{n+1}$ do not satisfy all conditions of Theorem
\ref{results:thm2}.  Thus $\xi_1$ must be algebraic over $K$ of
degree less than $n/t$.  Then, all of $\xi_1,\dots,\xi_{n+1}$ are
algebraic over $K$ of degree less than $n/t$, by virtue of the
recurrence relation which links these numbers.

%
%

\section{Construction of a polynomial}
 \label{constr}

 From this point on, the objective is to prove Proposition
\ref{mainprop:prop1}.  In this section, we fix a choice of real
numbers $X,Y\ge 1$ and assume that the convex body $\Cphi(X,Y)$
introduced in \S\ref{mainprop} contains a non-zero point
$\uy=(y_0,y_1,\dots,y_n)$ of $K^{n+1}$. We will derive several
consequences from this assumption.  Again, the constants
$c_9,c_{10},\dots$ that appear below, as well as implied constants
in the symbols $\ll$ and $\gg$, depend only on $K$, $n$, $w$ and
the points $\xi_1,\dots,\xi_t$.

For each integer $\ell=0,1,\dots,n$, we denote by $B_\ell\colon
E_\ell\times E_{n-\ell} \to K$ the $K$-bilinear form given by
\[
 B_\ell(F,G) = \varphi(FG,\uy)
\]
and we define $M_\ell$ to be the matrix of $B_\ell$ with respect
to the bases $\{1, T, \dots, T^\ell\}$ of $E_\ell$ and $\{1, T,
\dots, T^{n-\ell}\}$ of $E_{n-\ell}$.  Thus, $M_\ell$ is the
matrix of size $(\ell+1)\times (n-\ell+1)$ whose element of the
$i$-th row and $j$-th column is
\begin{equation*}
 B_\ell(T^{i-1},T^{j-1}) = \varphi(T^{i+j-2},\uy) = y_{i+j-2}
\end{equation*}
for $i=1,\dots,\ell+1$ and $j=1,\dots,n-\ell+1$.

Our first goal is to establish an upper bound for the height
$\rmH(M_\ell)$ of $M_\ell$ when $\ell\le n/2$, a condition which
ensures that $M_\ell$ has no more rows than columns (see
\S\ref{results} for the definition of the height). To this end, we
extend $B_\ell$ to a $K_w$-bilinear form $B_{\ell,w}\colon
E_{\ell,w}\times E_{n-\ell,w} \to K_w$ and we define $N_\ell$ to
be the matrix of $B_{\ell,w}$ with respect to the basis $\{1, T,
\dots, T^\ell\}$ of $E_{\ell,w}$ and the basis $\{R_0, R_1, \dots,
R_{n-\ell}\}$ of $E_{n-\ell,w}$ where
\begin{equation}
 \label{constr:R}
 R_0(T)=1
 \et
 R_j(T) = (T-\xi_1)(T-\xi_2)\cdots(T-\xi_j)
 \quad
 \text{for $j=1, \dots, n$,}
\end{equation}
extending for convenience the definition of $\xi_k$ for $k = t+1,
\dots, n$ by putting
\begin{equation}
 \label{constr:convention}
 \xi_{t+1}=\cdots=\xi_n=0.
\end{equation}

\begin{lemma}
 \label{constr:lemma1}
Let $\ell$ be an integer with $0\le \ell\le n/2$, and let $\uz_j$
denote the $j$-th column of $N_\ell$ for $j=1,\dots,n-\ell+1$.
Then, there are constants $c_9,c_{10},c_{11}\ge 1$ such that
\begin{itemize}
\item[(i)]
 $\|\uz_j\|_w \le c_9Y$ for $1\le j\le t$  and
 $\|\uz_j\|_w \le c_9X^{-1}$ for $t+1\le j\le n-\ell+1$,
\item[(ii)]
 $\rmH(M_\ell) \le c_{10} \|N_\ell\|_w$,
\item[(iii)]
 $\rmH(M_\ell) \le c_{11} Y^t X^{-(\ell+1-t)}$.
\end{itemize}
\end{lemma}

\begin{proof}
Fix an index $j$ with $1\le j\le n-\ell+1$.  By definition, we
have
\begin{equation*}
 \|\uz_j\|_w
 = \max_{0\le i\le \ell} |B_{\ell,w}(T^i,R_{j-1}(T))|_w
 = \max_{0\le i\le \ell} |\varphi_w(T^iR_{j-1}(T),\uy)|_w
 \le |\rho|_w^{-1}
\end{equation*} for any non-zero element $\rho$ of $K_w$ such that
$\rho T^iR_{j-1}(T)\in \cC_w(X,Y)$ for $i=0,\dots,\ell$.  As
$R_{j-1}(T)$ has bounded norm and as it is divisible by
$(T-\xi_1)\cdots(T-\xi_t)$ when $j> t$, there exists a constant
$c>0$ such that any choice of $\rho$ with $0<|\rho|_w \le cY^{-1}$
will do when $1\le j\le t$ and such that any choice of $\rho$ with
$0<|\rho|_w \le cX$ will do when $t+1\le j\le n-\ell+1$. The
inequalities (i) follow.

Now, fix a place $v$ of $K$ with $v\neq w$. Then, we have
$|\varphi_v(Q,\uy)|_v \le 1$ for any polynomial $Q\in E_{n,v}$
with $\|Q\|_v\le 1$. This implies that $\|\uy\|_v \le 1$ and
therefore that
\begin{equation*}
 \|M_\ell\|_v \le \max\{1,|(\ell+1)!|_v\}.
\end{equation*} Using the inequalities (i), we also find
\begin{equation*}
 \|N_\ell\|_w
 \le
 \max\{1,|(\ell+1)!|_w\} c_9^{\ell+1} Y^t X^{-(\ell+1-t)}
\end{equation*}
(this holds even when $\ell+1<t$). Since $M_\ell= N_\ell V$ for
some matrix $V\in\GL_{n-\ell+1}(K_w)$ with coefficients depending
only on $\xi_1,\dots,\xi_t$, $\ell$ and $n$, we also have $\|M_\ell\|_w
\le c' \|N_\ell\|_w$ for some constant $c'>0$.  Together with the
previous inequalities this proves (ii) with $c_{10} = n! c'$ and
(iii) with $c_{11} = n! c' c_9^n$.
\end{proof}

 To state the next result, we denote by $V_\ell$ the right kernel of
$B_\ell$ :
\begin{equation}
 \label{constr:Vell}
 V_\ell
 = \{ G \in E_{n-\ell} \,;\,
      \varphi(FG,\uy)=0 \text{ for all $F\in E_\ell$} \},
 \quad (0\le \ell \le n).
\end{equation}

\begin{lemma}
 \label{constr:lemma2}
Suppose that we have $c_{11}Y^t < X^{k+1-t}$ for some integer $k$
with $t\le k \le n/2$.  Then there exists an integer $h$ with
$1\le h\le k$ and a non-zero polynomial $P\in V_{n-h}$ which
divides any element of $V_{k-1}$ and satisfies
\begin{equation*}
 \deg(P) \le h
 \et
 \rmH(P)^{n-2h+2} \le c_{12} \rmH(M_{h-1})
\end{equation*} with a constant $c_{12}>0$ depending only on $n$.
\end{lemma}

\begin{proof}
By Lemma \ref{constr:lemma1} (iii), the condition $c_{11}Y^t <
X^{k+1-t}$ implies that $\rmH(M_k)<1$ and so $\rmH(M_k)=0$. Thus,
$M_k$ has rank at most $k$. Since $M_0$ has rank $1$, we deduce
that there exists an integer $h$ with $1\le h \le k$ such that
$\rank(M_{h-1})=h$ and $\rank(M_h)\le h$.  We now argue as in
Lemmas 6.3 and 6.4 of \cite{RW}.  Since $\rank(M_h)\le h$, the
left kernel of $B_h$ which is $V_{n-h}$ contains a non-zero
polynomial $P$. Then we have $\deg(P)\le h$ and $V_{h-1}$ contains
the set $P\cdot E_{n-2h+1}$ of all products $PQ$ with $Q\in
E_{n-2h+1}$. Since the dimension of $V_{h-1}$ is
$(n-h+2)-\rank(M_{h-1}) = n-2h+2$ and since $P\cdot E_{n-2h+1}$ is
a vector space over $K$ of the same dimension, we conclude that
$V_{h-1}=P\cdot E_{n-2h+1}$.  This equality has two consequences.
First of all, since $V_{k-1}$ is a subspace of $V_{h-1}$, the
polynomial $P$ divides all elements of $V_{k-1}$. Secondly, since
$\rmH(M_{h-1})$ is, by a well-known duality principle, equal to
the Schmidt height $\rmH(V_{h-1})$ of $V_{h-1}$, Proposition 5.2
of \cite{RW} shows that $\rmH(P)^{n-2h+2}\le c_{12} \rmH(M_{h-1})$
with a constant $c_{12}>0$ depending only on $n$ (see also
\S\ref{results}).
\end{proof}

We conclude this section by showing the following additional
property for the polynomial $P$ constructed in Lemma
\ref{constr:lemma2}:

\begin{lemma}
 \label{constr:lemma3}
Let the notation and the hypotheses be as in Lemma
\ref{constr:lemma2}.  Assume further that $k\le (n-t+2)/2$. Then
there exists an index $i$ with $1\le i\le t$ and an irreducible
factor $Q$ of $P$ such that
\begin{equation*}
 \left(
 \frac{|Q(\xi_i)|_w}{\|Q\|_w}
 \right)^t
 \le c_{13} X^{-\deg(Q)} \rmH(Q)^{-(n-2k+2)}.
\end{equation*}
\end{lemma}

\begin{proof}
Let $\uz_1, \uz_2,\dots,\uz_{n-h+2}$ denote the columns of $N_{h-1}$
and, for each $j=1,\dots,t+1$, let $N_{h-1}^{(j)}$ denote the
sub-matrix of $N_{h-1}$ whose columns are
$\uz_j,\uz_{j+1},\dots,\uz_{n-h+2}$. The hypothesis $k\le (n-t+2)/2$
ensures that each of these matrices has at least as many columns as
rows.  It also gives $h\le (n-t+2)/2$ which implies that the
products $T^iR_{j-1}(T)$ with $i=0,\dots,h-1$ and $j=1,\dots,t$ all
have degree at most $n-h$.  Since $P$ belongs to $V_{n-h}$, we
deduce that
\begin{equation}
 \label{constr:ortho}
 B_{h-1,w}(T^i,R_{j-1}(T)P(T))
 = B_{n-h,w}(T^iR_{j-1}(T),P(T))
 = 0
\end{equation}
for the same values of $i$ and $j$.

Fix an index $j$ with $1\le j\le t$.  Since $\deg(P)\le h$, there
exist a constant $c>0$ and elements $a_{j,1},\dots,a_{j,h}$ of
$K_w$ of absolute value at most $c\|P\|_w$ such that, with the
convention \eqref{constr:convention}, we have
\begin{equation*}
 P(\xi_j)-P(T)
 = \sum_{\ell=1}^h a_{j,\ell}(T-\xi_j)\cdots(T-\xi_{j+\ell-1}).
\end{equation*}
Taking into account that \eqref{constr:ortho} holds for $i=0,\dots
h-1$ we deduce that, for these values of $i$,
\begin{align*}
 P(\xi_j)B_{h-1,w}(T^i, R_{j-1}(T))
 &= B_{h-1,w}\big(T^i, R_{j-1}(T)(P(\xi_j)-P(T))\big) \\
 &= \sum_{\ell=1}^h a_{j,\ell} B_{h-1,w}(T^i,R_{j+\ell-1}(T)),
\end{align*}
and therefore that
\begin{equation*}
 P(\xi_j)\uz_j = \sum_{\ell=1}^h a_{j,\ell} \uz_{j+\ell}.
\end{equation*}
Applying this relation to all minors of order $h$ of
$N_{h-1}^{(j)}$ which include the column $\uz_j$ and using the
multilinearity of the determinant, this gives
\begin{align}
 |P(\xi_j)|_w \|N_{h-1}^{(j)}\|_w
 &\le \max\Big\{ |P(\xi_j)|_w, \sum_{\ell=1}^h |a_{j,\ell}|_w\Big\}
      \|N_{h-1}^{(j+1)}\|_w \\
 &\ll \|P\|_w \|N_{h-1}^{(j+1)}\|_w.
  \notag
\end{align}
Combining these relations for $j=1,\dots,t$, we get
\begin{equation*}
 \prod_{j=1}^t \frac{|P(\xi_j)|_w}{\|P\|_w}
 \le
 \frac{\|N_{h-1}^{(t+1)}\|_w}{\|N_{h-1}\|_w}.
\end{equation*}
On the other hand, since $\deg(P)\le h\le k\le n$, the estimates
of the two preceding lemmas give
\begin{equation*}
 \|N_{h-1}^{(t+1)}\|_w
 \ll \big(\max_{j>t} \|\uz_j\|_w\big)^h
 \ll X^{-h}
 \le X^{-\deg(P)}
\end{equation*}
and
\begin{equation*}
 \|N_{h-1}\|_w
 \gg \rmH(M_{h-1})
 \gg \rmH(P)^{n-2h+2}
 \ge \rmH(P)^{n-2k+2}.
\end{equation*}
So, if $i$ denotes an index for which $|P(\xi_i)|_w$ is minimal,
we find
\begin{equation*}
 \left( \frac{|P(\xi_i)|_w}{\|P\|_w} \right)^t
 \ll X^{-\deg(P)} \rmH(P)^{-(n-2k+2)}.
\end{equation*}
By multiplicativity, it follows that at least one irreducible
factor $Q$ of $P$ has the same property.
\end{proof}

%
%

\section{Degree and height estimates}
\label{est}

The notation being the same as in the previous section, our goal
is now to provide estimates for the degree and height of the
polynomial $Q$ constructed in Lemma \ref{constr:lemma3}.  We start
with the following construction of an auxiliary polynomial
(compare with \cite[Prop.\ 7.2]{RW}).

\begin{lemma}
 \label{est:lemma1}
Let the notation be as in Lemma \ref{constr:lemma3}. Then, there
exists a constant $c_{14}$ with $0<c_{14}<1$ such that the
following properties hold.
\par
{\rm (i)} Suppose that, for some integer $u\ge 0$, we have
\begin{equation}
 \label{est:condgen}
(XY)^{t+su} \le c_{14} X^{n-2k+3}.
\end{equation}
Then there exists a non-zero polynomial $G\in E_{n-2k+2}$ of
height at most $X$ such that $G^{(j)}$ belongs to $V_{k-1}$ for
$j=0,\dots,u$.
\par
{\rm (ii)} Suppose that, for some integer $u\ge 0$, we have
\begin{equation}
 \label{est:condpart}
(XY)^{t+u} \le c_{14} X^{n-2k+3},
\end{equation}
Suppose moreover that $\xi_1,\dots,\xi_t$ are all distinct, that
we have $\xi_1\notin K$ and $t\ge 2$, and that there exists a
polynomial $\linpol\in K[T]$ of degree $1$ such that
$\xi_{i+1}=\linpol(\xi_i)$ for $i=1,\dots,t-1$. Then there exists
a non-zero polynomial $G\in E_{n-2k+2}$ of height at most $X$ such
that $G\circ\linpol^j$ belongs to $V_{k-1}$ for $j=0,\dots,u$,
where $\linpol^j$ denotes the $j$-th iterate of $\linpol$.
\end{lemma}

\begin{proof}
The result follows from the adelic Minkowski convex body theorem
applied to an adelic convex body $\cK=\prod_v \cK_v$ of
$E_{n-2k+2}$ that we construct as follows, subject to the choice
of a real number $c$ with $0< c\le 1$.

In the case (i), we put $\cS = \cM_\infty$.  For each
$v\in\cM\setminus\{w\}$, we define $\cK_v$ to be the set of
elements $G$ of $E_{n-2k+2,v}$ with
\begin{equation}
 \label{est:S}
 \|G\|_v \le
 \begin{cases}
   1 &\text{if $v\notin \cS$,} \\
   c &\text{if $v\in \cS$,}
 \end{cases}
\end{equation}
and we define $\cK_w$ to be the set of polynomials $G\in
E_{n-2k+2,w}$ satisfying
\begin{equation*}
 \|G\|_w \le cX
 \et
 |G^{(j)}(\eta_i)|_w \le cY^{-1} \quad
 (1\le i\le s,\ 0\le j\le m_i+u-1).
\end{equation*}
We choose $c$ small enough, as a function of $n$ and $\max_{1\le
i\le t}|\xi_i|_w$, so that, for each $G\in\cK$, the products
$T^\ell G^{(j)}(T)$ with $0\le \ell \le k-1$ and $0\le j\le u$ all
belong to the interior of $\cC(X,Y)$. Then, for a suitable choice
of $c_{14}$, the condition \eqref{est:condgen} ensures that the
volume of $\cK$ is large enough so that, by Proposition
\ref{geom:minkowski}, $\cK$ contains a non-zero polynomial $G$ of
$E_{n-2k+2}$. For such a polynomial, we have $\rmH(G)\le X$ and,
by the remark following the proof of Proposition
\ref{geom:prop:duality}, we get
\begin{equation*}
 \varphi(T^\ell G^{(j)}(T),\uy) =0
 \quad
 (0\le \ell \le k-1,\ 0\le j\le u)
\end{equation*}
which shows that $G^{(j)}\in V_{k-1}$ for $j=0,\dots,u$.

In the case (ii), the hypotheses $\linpol(\xi_1)=\xi_2$ and
$\xi_1\notin K$ imply that $\linpol\in K+KT$ is uniquely
determined by $\xi_1$ and $\xi_2$.  We denote by $\cS$ the union
of $\cM_\infty$ with the finite set of places $v\in\cM$ for which
$\|\linpol\|_v > 1$. For each $v\in\cM\setminus\{w\}$, we define
$\cK_v$ to be the set of elements $G$ of $E_{n-2k+2,v}$ satisfying
\eqref{est:S}, and we define $\cK_w$ to be the set of all
polynomials $G\in E_{n-2k+2,w}$ satisfying
\begin{equation*}
 \|G\|_w \le cX
 \et
 |(G\circ\linpol^j)(\xi_1)|_w \le cY^{-1} \quad
 (0\le j\le t+u-1).
\end{equation*}
In this situation, we choose $c$ depending only on $n$,
$\max_{v\in\cM} \|A\|_v$ and $\max_{1\le i\le t}|\xi_i|_w$
so that, for each $G\in\cK$, the products $T^\ell
(G\circ \linpol^j)(T)$ with $0\le \ell \le k-1$ and $0\le j\le u$
all belong to the interior of $\cC(X,Y)$. Then, as in the
preceding case, a suitable choice of $c_{14}$ in the condition
\eqref{est:condpart} ensures that $\cK$ contains a non-zero
polynomial $G$ of $E_{n-2k+2}$, and any such polynomial $G$ has
the requested properties.
\end{proof}

We conclude with the following result.

\begin{lemma}
 \label{est:lemma2}
Suppose that the conditions (i) or (ii) of Lemma \ref{est:lemma1}
are fulfilled for some integer $u\ge 0$ and that, in the case
(ii), we have $\linpol^j(T)\neq T$ for $j=0,\dots,n-2k+2$. Then,
under the hypotheses of Lemma \ref{constr:lemma3}, the irreducible
polynomial $Q$ of $K[T]$ produced by Lemma \ref{constr:lemma3} has
\begin{equation}
 \label{est:est}
 \deg(Q) \le \frac{n-2k+2}{u+1}
 \et
 \rmH(Q) \le c_{15} X^{1/(u+1)}.
\end{equation}
\end{lemma}

\begin{proof}
Let $G$ be as in the conclusion of Lemma \ref{est:lemma1}.  Since
$Q$ is an irreducible factor of a polynomial $P$ which by Lemma
\ref{constr:lemma2} divides any element of $V_{k-1}$, the
polynomial $Q$ divides $G^{(j)}$ for $j=0,\dots,u$ in the case
(i), and $G\circ\linpol^j$ for $j=0,\dots,u$ in the case (ii). In
the case (i), we deduce that $Q^{u+1}$ divides $G$ and the
estimates \eqref{est:est} follow.

In the case (ii), the polynomials $Q, Q\circ\linpol^{-1}, \dots,
Q\circ\linpol^{-u}$ are irreducible factors of $G$ of the same
degree, and we have
\begin{equation*}
 c^{-j}\rmH(Q) \le \rmH(Q\circ\linpol^{-j}) \le c^j\rmH(Q)
 \quad
 (0\le j \le u),
\end{equation*}
for a constant $c\ge 1$ depending only on $n$ and $\rmH(\linpol)$.
Since the relation $\linpol(\xi_1)=\xi_2$ determines uniquely
$\linpol$, this constant $c$ ultimately depends only on $n$,
$\xi_1$ and $\xi_2$.   Let $m$ be the largest integer with $1\le
m\le u+1$ such that $Q, {Q\circ\linpol^{-1}}, \dots,
{Q\circ\linpol^{-(m-1)}}$ are two by two relatively prime.  Then,
the product $\prod_{j=0}^{m-1} Q\circ\linpol^{-j}$ divides $G$ and
we deduce that
\begin{equation}
 \label{est:estm}
 m\deg(Q)\le \deg(G)
 \et
 \rmH(Q)^m \ll \rmH(G).
\end{equation}
If $m=u+1$, this gives \eqref{est:est}. Assume now that $m\le u$.
Then $Q\circ\linpol^m$ is a multiple of $Q$ and therefore
$\linpol^m$ permutes the roots of $Q$.  In particular, there exist a
root $\alpha$ of $Q$ and an integer $i$ with $1\le i\le\deg(Q)$ such
that $\linpol^{mi}(\alpha)=\alpha$. Then, $Q$ divides
$\linpol^{mi}(T)-T$.  By hypothesis, the latter polynomial is
non-zero since by \eqref{est:estm} we have $mi\le n-2k+2$.  Hence,
$Q$ has degree $1$ and its height is equal to that of
$\linpol^{mi}(T)-T$.  This again gives \eqref{est:est} upon
noting that the condition \eqref{est:condpart} implies $u+1 \le
n-2k+2$ (since $X,Y\ge 1$ and $c_{14}<1$).
\end{proof}

%
%

\section{Proof of the main Proposition \ref{mainprop:prop1}}
 \label{proof}

Let the notation and hypotheses be as in the statement of Theorem
\ref{results:thm1} (resp.\ Theorem \ref{results:thm2}).  Define
accordingly $\nu=4Dst$ (resp.\ $\nu=4t$) as in the statement of
Proposition \ref{mainprop:prop1}.  Define also $k=[(n+2)t/\nu]$
where the brackets denote the integer part. Since $n\ge \nu\ge
4t$, this integer $k$ satisfies
\begin{equation*}
 t \le k \le \min\{n/2, (n-t+2)/2\}.
\end{equation*}
We also note that, in the situation of Theorem \ref{results:thm2},
the point $(\xi_1,\dots,\xi_t)$ is a zero of the prime ideal $\fp$
of $K[x_1,\dots,x_t]$ generated by the polynomials
$x_{i+1}-\linpol(x_i)$ for $i=1,\dots,t-1$.  To be consistent with
the notation of Theorem \ref{results:thm1}, we then put $D=1$,
since this ideal $\fp$ has degree $1$.  In both cases, we define
$u=Dt$. We also denote by $\xi$ the point of $\bP^t(\bC_w)$ with
homogeneous coordinates $\uxi = (1,\xi_1,\dots,\xi_t)$, and by
$\fP$ the homogeneous prime ideal of $K[x_0,\dots,x_t]$ which is
mapped to $\fp$ under the specialization $x_0\mapsto 1$.

Assume, by contradiction, that the conclusion of the proposition
does not hold and define $Y$ as a function of $X$ by
$Y=X^{(n+2-\nu)/\nu}$.  Then, there exists a real number $X_0\ge
1$ such that, for each $X\ge X_0$, the convex body $\Cphi(X,Y)$
defined in \S\ref{mainprop} contains a non-zero point $\uy$ of
$K^{n+1}$. Then, assuming that $X_0$ is sufficiently large, all
conditions of Lemmas \ref{constr:lemma2}, \ref{constr:lemma3},
\ref{est:lemma1} and \ref{est:lemma2} are fulfilled.  Indeed, for
$X$ sufficiently large, we find
\begin{equation*}
 c_{11} (XY)^t = c_{11} X^{(n+2)t/\nu} < X^{k+1},
\end{equation*}
and so the main condition $c_{11}Y^t < X^{k+1-t}$ of Lemma
\ref{constr:lemma2} is satisfied.  We also find
\begin{equation*}
 (XY)^{t+su} = X^{(n+2)t(1+Ds)/\nu} \le X^{(n+2)/2}
 \quad
 \Big( \text{\ resp.\ }
 (XY)^{t+u} = X^{(n+2)/2}\
 \Big),
\end{equation*}
while $n-2k+3 > (n+2)/2$.  So, the main condition
\eqref{est:condgen} (resp.\ \eqref{est:condpart}) of Lemma
\ref{est:lemma1} is satisfied for each sufficiently large $X$.
Thus, assuming $X_0$ sufficiently large, Lemma \ref{constr:lemma3}
provides, for each $X\ge X_0$, an irreducible polynomial $Q=Q_X$
of $K[T]$ and an index $i=i_X$ with $1\le i\le t$ such that
\begin{equation}
 \label{proof:Qxi}
 \left(\frac{|Q(\xi_i)|_w}{\|Q\|_w}\right)^t
 \le
 c_{13} X^{-\deg(Q)} \rmH(Q)^{-(n-2k+2)}.
\end{equation}
By Lemma \ref{est:lemma2}, this polynomial satisfies
\begin{equation*}
 \deg(Q) \le \frac{n-2k+2}{u+1}
 \et
 \rmH(Q) \le c_{15} X^{1/(u+1)}.
\end{equation*}
Moreover, since $[K(\xi_i):K] \ge n/u >\deg(Q)$, we also have
$Q(\xi_i) \neq 0$.  Define
\begin{equation*}
 n' = \left[\frac{n-2k+2}{u+1}\right]
 \et
 Y = c_{15} X^{1/(u+1)},
\end{equation*}
and let $P=P_Y$ denote the homogeneous polynomial of
$K[x_0,\dots,x_t]$ with the same degree as $Q$, for which
$P(1,x_1,\dots,x_t)=Q(x_i)$.  Then, $P$ has degree at most $n'$
and height at most $Y$.  It does not belong to $\fP$ since it does
not vanish at the point $\xi$.  Moreover, by \eqref{proof:Qxi}, we
have
\begin{equation*}
 \frac{ |P(\uxi)|_w }
      { \|P\|_w \|\uxi\|_w^{\deg(P)} }
 \ll
 X^{-\deg(P)/t} \rmH(P)^{-(n-2k+2)/t}
 \ll
 Y^{-1/t} \left( Y^{\deg(P)} \rmH(P)^{n'} \right)^{-D}
\end{equation*}
so that, for any sufficiently large value of $Y$, the above
polynomial $P_Y=P$ satisfies the main hypothesis
\eqref{gelfond:Pxi} of Theorem \ref{gelfond:thm1} below, with $n$
replaced by $n'$.  This is a contradiction as none of these
polynomials vanish at $\xi$.

%
%

\appendix
\section[A]{A version of Gel'fond's criterion for curves}
\label{gelfond}

In this appendix, we denote by $\bC_w$ the completion of $\Kbar$
with respect to its unique absolute value (also denoted $|\ |_w$)
which extends $|\ |_w$ on $K$. Then, $\bC_w$ is an algebraically
closed field containing $K_w$ as a subfield. We also fix a
positive integer $t$ and, for conciseness, we put $K[\ux] =
K[x_0,\dots,x_t]$ where $x_0,\dots,x_t$ denote independent
variables over $K$. We denote by $\deg(P)$ the degree of a
homogeneous polynomial $P$ of $K[\ux]$ and by $\rmH(P)$ its
height, that is the height of the vector of its coefficients.
Similarly, for a homogeneous ideal $I$ of $K[\ux]$, we denote by
$\deg(I)$ its degree and by $\rmH(I)$ the height of a Chow form of
$I$ (see below for a precise definition). Our goal is to prove the
following result which generalizes Theorem 4.2 of \cite{RW} (see
also Theorem 2b of \cite{DS}).

\begin{theorem}
 \label{gelfond:thm1}
Let $n$ be a positive integer, let $\fP$ be a homogeneous prime
ideal of $K[\ux]$ whose zero set $\Vbar$ in $\bP^t(\bC_w)$ has
dimension $1$, let $D=\deg(\fP)$, and let $\uxi = (\xi_0, \dots,
\xi_t)$ be homogeneous coordinates of a point $\xi$ of $\Vbar$.
Suppose that, for any sufficiently large real number $Y\ge 1$,
there exists a homogeneous polynomial $P=P_Y\in K[\ux]$ of degree
at most $n$ and height at most $Y$ which does not belong to $\fP$
and satisfies
\begin{equation}
 \label{gelfond:Pxi}
 \frac{|P(\uxi)|_w}{\|P\|_w\|\uxi\|_w^{\deg(P)}}
 <
 e^{-24 t^3 n^2 D} \rmH(\fP)^{-n^2}
 \big( \rmH(P)^n Y^{\deg(P)} \big)^{-D}.
\end{equation}
Then the point $\xi$ is defined over an algebraic extension of $K$
of degree at most $nD$ and the above polynomials vanish at this
point for any sufficiently large $Y$.
\end{theorem}

Our proof follows essentially the arguments of P.~Philippon in \S
II.3 of \cite{Ph}, taking advantage of a simpler context. For
convenience, we base this proof on the formalism and results of
Yu.~V.~Nesterenko in \cite{Nes}. We start by recalling the notion
of a Chow form of a homogeneous ideal $I$ of $K[\ux]$ and related
concepts which, in view of our present choice of normalization for
the absolute values of $K$ (see \S \ref{results}), differ slightly
from those of \cite{Nes}.

\medskip
Let $u_{i,j}$ for $i=1,\dots,t+1$ and $j=0,\dots,t$ be independent
variables over $K[\ux]$, and write $\uu_i=(u_{i,0},\dots,u_{i,t})$
for $i=1,\dots,t+1$.  Let $I$ be a homogeneous ideal of $K[\ux]$,
let $Z$ denote the set of zeros of $I$ in $\bP^t(\bC_w)$, and put
$r=\dim(Z)+1$ with the convention that $r=0$ if $Z$ is empty.
Denote by $I(r)$ the ideal of $K[\ux,\uu_1,\dots,\uu_r]$ generated
by the elements of $I$ and the polynomials $u_{i,0}x_0 + \cdots +
u_{i,t}x_t$ for $i=1,\dots,r$, and denote by $\Ibar(r)$ the ideal
of $K[\uu_1,\dots,\uu_r]$ consisting of the elements $G$ of that
ring for which there exists an integer $M\ge 1$ such that
$Gx_j^M\in I(r)$ for $j=0,\dots,t$.  Then $\Ibar(r)$ is a non-zero
principal ideal and we define a Chow form of $I$ to be any
generator $F$ of this ideal (see \S1 of \cite{Nes} for other
denominations and an historical perspective).  It is known that
such a polynomial of $K[\uu_1,\dots,\uu_r]$ is separately
homogeneous of degree $\deg(I)$ in each set of variables
$\uu_1,\dots,\uu_r$.  Moreover, as $F$ is uniquely determined up
to multiplication by a non-zero element of $K$, it makes sense (in
view of the product formula) to define the {\it height} $\rmH(I)$
of $I$ to be the height $\rmH(F)$ of the vector of coefficients of
$F$.

Let $S^{(1)},\dots,S^{(r)}$ be skew symmetric matrices of order
$t+1$ whose coefficients above the diagonal are altogether
independent variables over $\bC_w$, and let $\kappa$ denote the
$K$-linear ring homomorphism from $K[\uu_1, \dots, \uu_r]$ to
$\bC_w[S^{(1)}, \dots, S^{(r)}]$ mapping $\uu_i$ to $\uxi S^{(i)}$
for $i=1,\dots,r$. Following \cite{Nes}, we define the absolute
value of $I$ at $\xi$ by
\begin{equation*}
 |I(\xi)|_w
 =
 \frac{\|\kappa(F)\|_w}{\|F\|_w \|\uxi\|_w^{r\deg(I)}},
\end{equation*}
where $\|F\|_w$ (resp.\ $\|\kappa(F)\|_w$) stands for the largest
absolute value of the coefficients of F (resp.\ $\kappa(F)$). This
is independent of the choice of $F$ as well as the choice of
homogeneous coordinates $\uxi$ for $\xi$.  Moreover, we have
$|I(\xi)|_w=0$ if and only if $\xi$ belongs to an irreducible
component of $Z$ of dimension $r-1$.

Finally, we define the distance between $\xi$ and a point $z$ of
$\bP^t(\bC_w)$ with projective coordinates $\uz=(z_0,\dots,z_t)$
by the formula
\begin{equation*}
 \dist(\xi,z)
 =
 \|\uxi\|_w^{-1} \|\uz\|_w^{-1}
 \max_{0\le j,k\le t} |\xi_jz_k-\xi_kz_j|_w
\end{equation*}
(again this is independent of the choices of coordinates $\uxi$
for $\xi$ and $\uz$ for $z$).  Accordingly, we define the distance
between $\xi$ and the set $Z$ of zeros of $I$ in $\bP^t(\bC_w)$ by
\begin{equation*}
 \dist(\xi,Z) = \inf\{\dist(\xi,z)\,;\, z\in Z\}.
\end{equation*}

We can now state the results of \cite{Nes} that we need.

\begin{lemma}
 \label{gelfond:lemma1}
Let $J$ be an unmixed homogeneous ideal of $K[x_0,\dots,x_t]$ and
let $\fp_1,\dots,\fp_s$ be its associated prime ideals. Then,
there exist integers $k_1,\dots,k_s\ge 1$ such that
\begin{itemize}
\item[(i)]   $\sum_{j=1}^s k_j \deg(\fp_j) \le \deg(J)$,
\item[(ii)]  $\prod_{j=1}^s \rmH(\fp_j)^{k_j}
              \le e^{t^2\deg(J)} \rmH(J)$,
\item[(iii)] $\prod_{j=1}^s |\fp_j(\uxi)|_w^{k_j}
              \le e^{t^3\deg(J)} |J(\uxi)|_w$.
\end{itemize}
\end{lemma}

This follows immediately from a simple adaptation of the proof of
Proposition 1.2 of \cite{Nes} upon noting that, for each
Archimedean place $v\in\cM_\infty$, the absolute value $|{\ }|_v^{
d/d_v }$ of $K$ coincides with the usual absolute Archimedean
value on $\bQ$ and that we have $\sum_{v\in\cM_\infty} d_v/d =1$.
Then, the local estimates of \cite{Nes} applied to the absolute
values $|{\ }|_v^{d/d_v}$ with $v\in\cM_\infty$ combine to give
assertion (ii) and lead to a stronger form of the assertion (iii)
where the argument of the exponential is replaced by $(d_w/d) t^3
\deg(J)$.

\begin{lemma}
 \label{gelfond:lemma2}
Let $\fq$ be a homogeneous prime ideal of $K[\ux]$ whose zero set
$Z$ in $\bP^t(\bC_w)$ is not empty, and let $P$ be a homogeneous
polynomial from $K[\ux]$ with $P\notin \fq$. Put
\begin{equation*}
 r=\dim(Z)+1, \quad
 \rho=\dist(\xi,Z)
 \et
 \delta= \frac{|P(\uxi)|_w}{\|P\|_w\|\uxi\|_w^{\deg(P)}}.
\end{equation*}
\par%
If $r\ge 2$, there exists a homogeneous unmixed ideal $J$ of
$K[\ux]$ with the following properties.  Its set of zeros in
$\bP^t(\bC_w)$ has dimension $r-2$ and coincide with that of
$(\fq,P)$. Moreover, we have:
\begin{itemize}
\item[(i)]
 $\deg(J) \le \deg(\fq)\deg(P)$,
\item[(ii)]
 $\rmH(J) \le e^{2t^2\deg(\fq)\deg(P)}
              \rmH(\fq)^{\deg(P)} \rmH(P)^{\deg(\fq)}$,
\item[(iii)]
 $|J(\uxi)|_w \rmH(J)
  \le
  e^{11t^2\deg(\fq)\deg(P)}
  \rmH(\fq)^{\deg(P)} \rmH(P)^{\deg(\fq)}
  \begin{cases}
    \delta
    &\mbox{if }\rho < \delta,\\
    |\fq(\uxi)|_w
    &\mbox{otherwise.}
  \end{cases}$
\end{itemize}
\par%
If $r=1$, the above inequality (iii) holds with the left hand side
replaced by $1$.
\end{lemma}

This second result follows from a similar adaptation of the proof
of Proposition 1.4 of \cite{Nes}.  Part (ii) uses moreover $r+1\le
2t$ while part (iii) requires replacing the inequality (37) on
page 314 of \cite{Nes} with an equality involving the height of
the given Chow form $G$ of $J$.

\begin{corollary}
 \label{gelfond:corollary}
Let $\fq$ and $P$ be as in Lemma \ref{gelfond:lemma2}.  Then, any
minimal prime ideal $\fp$ of $(\fq,P)$ satisfies
\begin{equation*}
 \deg(\fp) \le \deg(\fq)\deg(P)
 \et
 \rmH(\fp) \le e^{3t^2\deg(\fq)\deg(P)}
               \rmH(\fq)^{\deg(P)} \rmH(P)^{\deg(\fq)}.
\end{equation*}
\end{corollary}

\begin{proof}
Lemma \ref{gelfond:lemma2} provides a homogeneous unmixed ideal
$J$ of $K[\ux]$ whose minimal prime ideals are the same as those
of $(\fq,P)$. Since the degree of any ideal is bounded below by
$0$ while its height is bounded below by $1$, Lemma
\ref{gelfond:lemma1} then shows that each minimal prime ideal
$\fp$ of $(\fq,P)$ satisfies $\deg(\fp)\le \deg(J)$ and
$\rmH(\fp)\le e^{t^2\deg(J)} \rmH(J)$. The conclusion follows
using the upper bounds for $\deg(J)$ and $\rmH(J)$ provided by
Lemma \ref{gelfond:lemma2} (i) and (ii).
\end{proof}

\begin{proof}[Proof of Theorem \ref{gelfond:thm1}]
Choose $Y_0\ge 1$ such that $P_Y$ is defined for each $Y\ge Y_0$.
Then fix an arbitrary choice of $Y$ with $Y\ge Y_0$ and put
$P=P_Y$. Since $P\notin \fP$ and since $\dist(\xi,\Vbar)=0$, Lemma
\ref{gelfond:lemma2} provides us with a homogeneous unmixed ideal
$J$ of $K[\ux]$ whose set of zeros in $\bP^t(\bC_w)$ has dimension
$0$ and coincide with that of $(\fP,P)$.  The same lemma also
gives estimates for this ideal, which taking into account the
inequality \eqref{gelfond:Pxi} and the fact that $\deg(P)\le n$,
imply $\deg(J)\le D\deg(P) \le nD$ and
\begin{align*}
  |J(\uxi)|_w \rmH(J)^n
  &\le
  e^{11 t^2 n^2 D} \rmH(\fP)^{n^2} \rmH(P)^{nD}
    \frac{|P(\uxi)|_w}{\|P\|_w\|\uxi\|_w^{\deg(P)}} \\
  &<
  e^{-13 t^3 n^2 D} Y^{-D\deg(P)}.
\end{align*}
Let $\fp_1,\dots,\fp_s$ be the associated prime ideals of $J$ in
$K[\ux]$. According to lemma \ref{gelfond:lemma1}, there exist
integers $k_1,\dots,k_s\ge 1$ such that $\sum_{j=1}^s k_j
\deg(\fp_j) \le \deg(J)$ and
\begin{equation*}
 \prod_{j=1}^s \big(|\fp_j(\uxi)|_w \rmH(\fp_j)^n\big)^{k_j}
 \le
  e^{t^2(t+n)\deg(J)} |J(\uxi)|_w \rmH(J)^n.
\end{equation*}
Putting these estimates together, we get
\begin{equation*}
 \prod_{j=1}^s \big(|\fp_j(\uxi)|_w \rmH(\fp_j)^n\big)^{k_j}
 <
 \big( e^{-11t^2n}Y^{-1} \big)^{\deg(J)}
 \le
 \prod_{j=1}^s \big( e^{-11t^2n}Y^{-1} \big)^{k_j \deg(\fp_j)}.
\end{equation*}
Therefore, there is at least one index $j$ for which the prime
ideal $\fp=\fp_j$ satisfies
\begin{equation}
 \label{gelfond:dist}
 |\fp(\uxi)|_w < e^{-11t^2nN} \rmH(\fp)^{-n} Y^{-N}
 \quad
\hbox{where $N=\deg(\fp)$.}
\end{equation}

Now, define $X_0$ as the infimum of all real numbers $X$ with
$X\ge Y_0$ such that $P_X\in \fp$.  By construction, we have $Y_0
\le X_0 \le Y$.  Moreover, for each $X\ge X_0$ such that $P_X\in
\fp$, the ideal $\fp$ is a minimal prime ideal of $(\fP,P_X)$.
Therefore, Corollary \ref{gelfond:corollary} gives
\begin{equation*}
 \rmH(\fp) \le e^{3t^2D\deg(P_X)} \rmH(\fP)^{\deg(P_X)} \rmH(P_X)^D
        \le e^{3t^2nD} \rmH(\fP)^n X^D.
\end{equation*}
As $X$ can be taken arbitrarily close to $X_0$, this implies
\begin{equation}
 \label{gelfond:height}
 \rmH(\fp) \le e^{3t^2nD} \rmH(\fP)^n X_0^D.
\end{equation}

Assume for the moment that $Y_0 < X_0$.  Choose $X$ with
$\max\{Y_0,X_0/2\} \le X < X_0$ and put $Q=P_X$. Define also
\begin{equation*}
 \rho= \dist(\xi,Z)
 \et
 \delta=\frac{|Q(\uxi)|_w}{\|Q\|_w \|\uxi\|^{\deg(Q)}}
\end{equation*}
where $Z$ denotes the zero set of $\fp$ in $\bP^t(\bC_w)$. Since
$Q\notin \fp$ and $\dim(Z)=0$, the inequality of Lemma
\ref{gelfond:lemma2} (iii) applies with the left hand side
replaced by $1$, $\fq$ replaced by $\fp$, and $P$ replaced by $Q$.
If $\rho<\delta$, then, taking into account \eqref{gelfond:height}
together with $N\le nD$, $\deg(Q)\le n$ and $X_0\le 2X$, this
gives
\begin{equation*}
 1 \le e^{11t^2nN} \rmH(\fp)^{\deg(Q)} \rmH(Q)^N \delta
   \le e^{15 t^2 n^2 D} \rmH(\fP)^{n^2} X^{D\deg(Q)}
         \rmH(Q)^{nD} \delta
\end{equation*}
against the upper bound for $\delta$ associated with $Q=P_X$.  So,
we have $\rho\ge \delta$ and Lemma \ref{gelfond:lemma2} (iii) then
gives
\begin{equation*}
 1 \le e^{11t^2nN} \rmH(\fp)^{\deg(Q)} \rmH(Q)^N |\fp(\uxi)|_w.
\end{equation*}
which now contradicts \eqref{gelfond:dist} since $\deg(Q)\le n$
and $\rmH(Q)\le X\le Y$.

Thus, we have $X_0=Y_0$ which in view of \eqref{gelfond:height}
means that the height of $\fp$ is bounded above by a constant
which is independent of $Y$.  Since the degree of $\fp$ is bounded
by $nD$, this implies that, as $Y$ varies, $\fp$ stays within a
finite set of ideals.  Since the upper bound for $|\fp(\uxi)|_w$
given by \eqref{gelfond:dist} tends to zero as $Y$ tends to
infinity, this shows that $|\fp(\uxi)|_w=0$ for any sufficiently
large $Y$ and thus that $\xi$ is a zero of $\fp$ for those values
of $Y$.  In particular, we deduce that $\xi$ is defined over an
algebraic extension of $K$ of degree at most $nD$ and that
$P_Y(\xi)=0$ for any sufficiently large $Y$.
\end{proof}


\end{document}